\documentclass{article}
\usepackage[notref,notcite]{showkeys}
\usepackage{color}
\usepackage{amsmath,amssymb}
\usepackage[titletoc]{appendix}
\catcode`\@=11 \@addtoreset{equation}{section}
\def\theequation{\thesection.\arabic{equation}}
\catcode`\@=12

\newtheorem{Theorem}{Theorem}[section]
\newtheorem{Proposition}{Proposition}[section]
\newtheorem{Lemma}{Lemma}[section]
\newtheorem{Remark}{Remark}[section]
\newtheorem{Corollary}{Corollary}[section]
\newtheorem{Definition}{Definition}[section]

\newcommand{\bTheorem}[1]{
\begin{Theorem} \label{T#1} }
\newcommand{\eT}{\end{Theorem}}

\newcommand{\bProposition}[1]{
\begin{Proposition} \label{P#1}}
\newcommand{\eP}{\end{Proposition}}

\newcommand{\bLemma}[1]{
\begin{Lemma} \label{L#1} }
\newcommand{\eL}{\end{Lemma}}

\newcommand{\bCorollary}[1]{
\begin{Corollary} \label{C#1} }
\newcommand{\eC}{\end{Corollary}}

\newcommand{\bFormula}[1]{
\begin{equation} \label{#1}}
\newcommand{\eF}{\end{equation}}

\newcommand{\Om}{\Omega}

\newcommand{\DC}{C^\infty_c}

\newcommand{\vr}{\varrho}

\newcommand{\vrt}{\tilde{\vr}}

\newcommand{\vut}{\tilde{\vu}}
\newcommand{\ut}{\tilde{u}}

\newcommand{\vu}{\vc{u}}
\newcommand{\vv}{\vc{v}}
\newcommand{\vz}{\vc{z}}
\newcommand{\vU}{\vc{U}}
\newcommand{\vw}{\vc{w}}

\newcommand{\vF}{\vc{F}}
\newcommand{\vV}{\vc{V}}
\newcommand{\ve}{\vc{e}}
\newcommand{\vc}[1]{{\bf #1}}
\newcommand{\qed}{\rightline{ $\square$}}
\renewcommand{\div}{{\rm div}}
\newcommand{\Div}{{\rm div}_x}
\newcommand{\Grad}{\nabla_x}

\newcommand{\tn}[1]{\mbox {\F #1}}
\newcommand{\dx}{{\rm d} {x}}
\newcommand{\dS}{{\rm d} {S}}
\newcommand{\vph}{{\boldsymbol \varphi}}

\renewcommand{\d}{{\rm d} }
\newcommand{\dt}{{\rm d} t }

\newcommand{\dxdt}{\dx\dt}

\newcommand{\bProof}{{\bf Proof: }}

\newcommand{\n}{\vc{n}}

\newcommand{\vX}{\vc{X}}

\newcommand{\Ecal}{{\cal E}}
\newcommand{\Rcal}{{\cal R}}
\newcommand{\ep}{\varepsilon}
\newcommand{\de}{\partial}
\newcommand{\vn}{\vc{n}}
\newcommand{\ddiv}{{\rm div}\,}

\font\F=msbm10 scaled 1000
\newcommand{\R}{\mbox{\F R}}

\providecommand{\abs}[1]{\left\lvert#1\right\rvert}
\newcommand{\sil}{\rightarrow}
\newcommand{\norm}[1]{\ensuremath{\left\| #1 \right\|}}
\newcommand{\norma}[2]{\ensuremath{\norm{#1}_{#2}}}

\definecolor{grey}{rgb}{0.85,0.85,0.85}
\date{}

\long\def\greybox#1{%
    \newbox\contentbox%
    \newbox\bkgdbox%
    \setbox\contentbox\hbox to \hsize{%
        \vtop{
            \kern\columnsep
            \hbox to \hsize{%
                \kern\columnsep%
                \advance\hsize by -2\columnsep%
                \setlength{\textwidth}{\hsize}%
                \vbox{
                    \parskip=\baselineskip
                    \parindent=0bp
                    #1
                }%
                \kern\columnsep%
            }%
            \kern\columnsep%
        }%
    }%
    \setbox\bkgdbox\vbox{
        \color{grey}
        \hrule width  \wd\contentbox %
               height \ht\contentbox %
               depth  \dp\contentbox
        \color{black}
    }%
    \wd\bkgdbox=0bp%
    \vbox{\hbox to \hsize{\box\bkgdbox\box\contentbox}}%
    \vskip\baselineskip%
}

\begin{document}


\title{Local existence of strong solutions and weak-strong uniqueness for the compressible Navier-Stokes system on moving domains}

\author{
Ond{\v r}ej Kreml $^1$\thanks{The  work  of O.K.   was partially supported by the project 7AMB16PL060 of Czech Ministry of Education, Youth and Sports, in the framework of RVO:67985840}
\and {\v S}{\' a}rka Ne{\v c}asov{\' a} $^1$\thanks{The work of \v S.N.  was partially supported by the project 7AMB16PL060 of Czech Ministry of Education, Youth and Sports(last section of paper). The rest  of paper was  supported by  the Czech Science Foundation project  GA17-01747 in the framework of RVO:67985840}\footnotemark[1]
\and Tomasz Piasecki $^2$\thanks{The work of T.P was partially supported by the polish National Science Centre Harmonia project UMO-2014/14/M/ST1/00108 and his stay in Institute of Mathematics by project 7AMB16PL060.}
}

\maketitle

\bigskip

\centerline{$^1$Institute of Mathematics, Czech Academy of Sciences} \centerline{\v Zitn\' a 25, 115 67 Praha,
Czech Republic}
\bigskip

\centerline{$^2$  Institute of Applied Mathematics and Mechanics, University of Warsaw}
\centerline{Banacha 2, 02-097 Warszawa, Poland}

\medskip

\begin{abstract}

We consider the compressible Navier-Stokes system on
time-dependent domains with prescribed motion of the boundary. For both the no-slip boundary conditions as well as slip boundary conditions
we prove local-in-time existence of strong solutions. These results are obtained using a transformation
of the problem to a fixed domain and an existence theorem for Navier-Stokes like systems with lower order terms
and perturbed boundary conditions. We also show the weak-strong uniqueness principle for slip boundary conditions which remained so far open question.

\end{abstract}

\medskip

\noindent
{\bf Keywords:} compressible Navier-Stokes equations, time-dependent domain, strong solution, local existence, weak-strong uniqueness.

\medskip

\noindent
{\bf MSC:} 35Q30, 76N10

\section{Introduction}\label{i}

We consider a barotropic flow of a compressible viscous fluid in the absence of external forces. Such flow is described by the isentropic compressible Navier-Stokes system
\bFormula{i1a}
\partial_t \vr + \Div (\vr \vu) = 0,
\eF
\bFormula{i1b}
\partial_t (\vr \vu) + \Div (\vr \vu \otimes \vu) + \Grad p(\vr) = \Div \tn{S}(\Grad \vu),
\eF
where $\vr$ is the density of the fluid and $\vu$ denotes the velocity. We assume that the stress tensor $\tn{S}$ is determined by the standard Newton rheological law
\bFormula{i4}
\tn{S} (\Grad \vu) = \mu \left( \Grad \vu +
\Grad^t \vu - \frac{2}{3} \Div \vu \tn{I} \right) + \eta \Div \vu \tn{I}
\eF
with $\mu > 0$ and $\eta \geq 0$ being constants. The pressure $p(\vr)$ is a given sufficiently smooth function of the density and
we introduce the pressure potential as
\bFormula{p1pb}
H(\vr) = \vr\int_1^\vr \frac{p(z)}{z^2} {\rm d} z.
\eF
We are interested in proving local existence of strong solutions to the system of equations \eqref{i1a}-\eqref{i1b} on a moving domain $\Omega = \Omega_t$ with prescribed movement of the boundary. More precisely, the boundary of the domain $\Omega_t$ is assumed to be described by means of a given velocity field $\vc{V}(t,x)$, where $t \geq 0$ and
$x \in \mathbb{R}^3$. Assuming $\vc{V}$ is regular and solving the associated system of differential equations
\begin{equation}\label{eq:coc}
\frac{{\rm d}}{{\rm d}t} \vc{X}(t, x) = \vc{V} \Big( t, \vc{X}(t, x) \Big),\ t > 0,\ \vc{X}(0, x) = x
\end{equation}
we set
\[
\Omega_\tau = \vc{X} \left( \tau, \Omega_0 \right),\]
where $\Omega_0 \subset \mathbb{R}^3$ is a given domain. Moreover we denote $\Gamma_\tau = \partial \Omega_\tau$ and
\[ Q_\tau = \bigcup_{t \in (0,\tau)} \{t\}\times \Omega_t =: (0,\tau)\times \Omega_t.
\]

System \eqref{i1a}-\eqref{i1b} is supplied with the Navier slip boundary condition stated as the combination of the impermeability relation
\bFormula{i6b} (\vu - \vc{V})
\cdot \vc{n} |_{\Gamma_\tau} = 0 \ \mbox{for any}\ \tau \geq 0,
\eF
where $\vc{n}(\tau,x)$ denotes the unit outer normal vector to the boundary $\Gamma_\tau$, and equation
\bFormula{i6c}
\left(\left[ \tn{S} \vc{n} \right]_{\rm tan} + \kappa\left[ \vu - \vc{V} \right]_{\rm tan}\right)|_{\Gamma_\tau} = 0,
\eF
where $\kappa \geq 0$ represents a friction coefficient. In particular the choice $\kappa = 0$ yields the complete slip boundary condition and in the limit $\kappa \rightarrow +\infty$ we recover the no-slip boundary condition \eqref{i6a}.
Notice that denoting $\tau_k$, $k= 1, 2$, the tangent vectors to the boundary $\Gamma_\tau$, we have $\Div \vu\, \tn{I}\,\vn\cdot \tau_k = 0$ and therefore condition \eqref{i6c} reduces to 
\begin{equation}
\left(\vn \cdot \mu(\Grad \vu +\Grad^t \vu)\cdot \tau_k+\kappa (\vu-\vV)\cdot \tau_k\right)|_{\Gamma_\tau} = 0, \quad k=1,2
\end{equation}
%

Finally, the system of equations \eqref{i1a}-\eqref{i1b} is supplemented by the initial conditions
\bFormula{i7}
\vr(0, \cdot) = \vr_0 ,\quad
\vu(0, \cdot) = \vu_0 \quad \mbox{in}\ \Omega_0.
\eF

Global weak solutions to the compressible barotropic Navier--Stokes system on a fixed domain were proved 
to exist in a pioneering work by Lions \cite{LI4}. This theory was later extended by Feireisl and collaborators 
(\cite{FNP}, \cite{EF61}, \cite{EF70}, \cite{EF71}) to cover larger class of pressure laws. 
The existence theory in the case of moving domains was developed in \cite{FeNeSt} 
for no-slip boundary conditions using the so-called Brinkman penalization and in \cite{FKNNS} 
for slip boundary conditions. Recently these results have been generalized to a complete system with 
thermal effects in \cite{KMNW1} and \cite{KMNW2}.

There is also quite large amount of literature available concerning existence of strong solutions for the system \eqref{i1a}--\eqref{i1b} 
(or even for more complex systems involving also heat conductivity of the fluid) on a fixed domain, 
such solutions are proved to exist either locally in time or globally provided the initial data are 
sufficiently close to a rest state, let us mention among others 
\cite{MN1}, \cite{MN2}, \cite{V1}, \cite{V2} in the framework of Hilbert spaces and \cite{MZ1},\cite{MZ2},\cite{MZ3} in the $L^p$ setting. All these results are proved under the no-slip boundary conditions.
The case of slip boundary condition on fixed domain was considered by Zaj\c aczkowski \cite{Za} and Hoff \cite{Ho}. 
Free boundary problems for the system \eqref{i1a}-\eqref{i1b} has been investigated by Zaj\c aczkowski et al. 
(\cite{Za1},\cite{Za2}) where global existence of strong solutions in $L^2$ setting is shown under assumption that the domain is close to a ball. 
Enomoto and Shibata \cite{ES} have shown local existence and global existence for small data using entirely different approach in the maximal $L_p-L_q$ regularity setting for Dirichlet boudary conditions. Analogous results for slip bounadry conditions have been shown in \cite{Mur} and \cite{SM}.
For moving domains with given motion of the boundary local-in-time existence results of strong solutions (incompressible case) can be found in the $L^p$ setting by Hieber and al. \cite{DGH}. Moreover, in the case of fluid-structure interaction we can mention work of Hieber for incompressible and also compressible case \cite{GGH,HM}.

The concept of relative entropies has been successfully used in the context of partial differential equations 
(see among others Carrillo et al. \cite{Ca}, Masmoudi \cite{MA}, Saint-Raymond \cite{SR}, 
Wang and Jiang \cite{WaJa}). Germain \cite{Ge} introduced a notion of solution to the 
system \eqref{i1a}--\eqref{i1b} based on a relative entropy inequality with respect to a hypothetical 
strong solution. Similar idea was adapted by Feireisl et al. \cite{FeNoSu} who defined a suitable weak 
solution to the barotropic Navier-Stokes system based on a general relative entropy inequality with 
respect to any sufficiently smooth pair of functions. In \cite{FeJiNo} the authors used relative entropy 
inequality to prove the weak-strong uniqueness property.
In the case of moving domain with no-slip boundary condition   Doboszczak \cite{dobo} proved both 
the relative entropy inequality as well as the weak-strong uniqueness property under assumption 
of an existence of local strong solution.

Aim of our paper is to fill the gap in the theory of weak-strong uniqueness. First we prove the local existence 
of strong solutions to the system \eqref{i1a}--\eqref{i1b}. The general approach consist in nowadays classical 
method of splitting the problem into continuity and momentum equations and investigation of linearized problems. 
Here we adapt {\it a mixed approach}. The reason is that standard approach to change coordinate system to 
fixed domain is not possible since we cannot close our fixed point argument. 
It enforces us to solve the continuity equation directly on moving domain. 
For this purpose we extend the known existence theory for the linear transport system  equation based on the 
method of characteristics. Analogous result on a fixed domain has been applied for example in \cite{Za}, 
however without proof which turns out to be quite involved and here we present it for the sake of completeness,
already on a moving domain (Proposition 3.1). To treat the linear momentum equation we apply the 
Lagrangian coordinates. The main difficulty in the resulting system on a fixed domain is in nonhomogeneous 
boundary conditions which, in spite of smallness of time, require careful treatment with appropriate extension 
of the boundary data. The result (Lemma 4.2) is based on explicit construction and careful 
estimates. We shall emphasize that the proof does not rely on particular structure of our boundary data, 
therefore it can be adapted to other nonhomogeneous boundary value problems.
Having completed the proof of the local existence we show that weak solutions to the system \eqref{i1a}--\eqref{i1b} on moving domains proved to exist in \cite{FKNNS} possess also an energy inequality (Proposition \ref{p:ex2}). Using this energy inequality we prove the relative entropy (energy) inequality (Proposition \ref{p:REI}) and finally the weak-strong uniqueness property (Theorem \ref{t:WSU}). 
Novelty in our paper is that we use mixed approach in the proof of local existence, 
developing a direct approach to linear transport equation on moving domains, 
and the construction of extension of the boundary data.
Both results can be of independent interest and find applications in other problems.
Moreover, we complete the theory of  weak-strong uniqueness on moving domain with Navier type of boundary 
condition.

We are now in a position to formulate the main results of our paper.
The first one gives the local existence of strong solutions (a precise definition of function spaces defined on the moving domain is given in Section \ref{ss:spaces}).
\bTheorem{main}
Let $\Omega_0 \subset \mathbb{R}^3$ be a bounded domain of class $C^2$.
Assume that $p(\vr)$ is a $C^2$ function of the density
and $\vc{V} \in C^3((0,T) \times \mathbb{R}^3)$.
Assume further that $\vu_0 \in H^3(\Omega_0), \vr_0 \in H^2(\Omega_0)$ and there exists positive constants $c_1,c_2$ such that
$0<c_1 \leq \vr_0 \leq c_2$. Then there exists (sufficiently small) $T>0$ and a
unique solution $(\vu,\vr)$ to the system \eqref{i1a}-\eqref{i4}
with boundary conditions \eqref{i6b}-\eqref{i6c} (or \eqref{i6a}) and initial condition \eqref{i7}
such that\\
$\vu \in L^\infty(0,T;H^2(\Omega_t))\cap L^2(0,T;H^3(\Omega_t)), \vu_t \in L^\infty(0,T;H^1(\Omega_t)) \cap L^2(0,T;H^2(\Omega_t)),
\vu_{tt} \in L^2(0,T;L^2(\Omega_t))$
$\vr \in L^\infty(0,T;H^2(\Omega_t)), \vr_t \in L^2(0,T;H^1(\Omega_t))$.
\eT
The second main result concerns the weak-strong uniqueness principle.
\begin{Theorem}\label{t:WSU}
let $\vc{V} \in C^1([0,T]; C^{3}_c (\mathbb{R}^3))$ be given.
Assume that the pressure $p \in C[0, \infty) \cap C^1(0, \infty)$ satisfies
\[
p(0) = 0,\ p'(\vr) > 0 \ \mbox{for any}\ \vr > 0,\ \lim_{\vr \to \infty} \frac{p'(\vr)}{\vr^{\gamma - 1}} = p_\infty > 0
\ \mbox{for a certain}\ \gamma > 3/2.
\]
Let $(\vr,\vu)$ be a weak solution to the compressible Navier-Stokes system \eqref{i1a}-\eqref{i7} constructed in Theorem \ref{t:ex}. Let $(\vrt,\vut)$ be a strong solution to the same problem satisfying
\begin{equation}
0 < \inf_{Q_T} \vrt \leq \sup_{Q_T} \vrt < \infty
\end{equation}
\begin{equation} \label{Lp}
\Grad \vrt \in L^2(0,T;L^q(\Omega_t)), \quad \Grad^2 \vut \in L^2(0,T;L^q(\Omega_t))
\end{equation}
with $q > \max\{3;\frac{6\gamma}{5\gamma-6}\}$, and emanating from the same initial data. Then, up to a set of measure zero,
\begin{equation}
\vr = \vrt, \, \vu = \vut \quad \text{ in } Q_T.
\end{equation}
\end{Theorem}
\begin{Remark}
Notice that the strong solution constructed in Theorem \ref{Tmain} satisfies the assumptions of Theorem \ref{t:WSU} since by the imbedding theorem we have \eqref{Lp} for $q \leq 6$ and $\frac{6\gamma}{5\gamma-6}<6$ for $\gamma>\frac{3}{2}$.
\end{Remark}
The paper is organized as follows. In Section \ref{s:prel} we discuss the proper definition of function spaces on moving domains, recall the existence theorem for weak solutions from \cite{FKNNS} and introduce the iterative scheme used in the proof of Theorem \ref{Tmain}. In Sections \ref{s:LinCont} and \ref{s:LinMom} we present the existence theory for the linear continuity and momentum equations respectively in appropriate function spaces. In Section \ref{s:proof} we show the convergence of the iterative scheme which completes the proof of Theorem \ref{Tmain}. Section \ref{s:WSU} is dedicated to the proof of Theorem \ref{t:WSU}. Finally, in Section \ref{s:conc} we present concluding remarks regarding different boundary conditions, regularity of $\vV$ etc.

\section{Preliminaries}\label{s:prel}

\subsection{Function spaces}\label{ss:spaces}

To begin with  we introduce the function spaces $L^p(0,T,X(\Omega_t))$. For fixed $T > 0$ we assume that there exists $R > 0$ such that for all $t \in [0,T]$ it holds $\Omega_t \subset B_R(0)$, where $B_R(0)$ denotes the ball in $\mathbb{R}^3$ of radius $R$ centered at the origin. Then we define
\begin{equation}\label{prostor}
L^p(0,T;L^q(\Omega_t)) := \left\{u \in L^p(0,T,L^q(B_R(0))), u(t,\cdot) = 0 \text{ in } B_R(0) \setminus \Omega_t \; \text{for a.e. } t\in(0,T)\right\}
\end{equation}
with the norm
\[
\norma{u}{L^p(0,T,L^q(\Omega_t))} := \left(\int_0^T \norma{u(t)}{L^q(\Omega_t)}^p dt\right)^{\frac 1p}
\]
for $p < \infty$ and
\[
\norma{u}{L^\infty(0,T,L^q(\Omega_t))} := \text{ess\, sup}_{t \in (0,T)} \norma{u(t)}{L^q(\Omega_t)}.
\]

Similarly we define spaces $L^p(0,T;W^{l,q}(\Omega_t))$. Let $l \in \mathbb{N}$ and $\alpha$ be a multi-index. Then
\[
L^p(0,T,W^{l,q}(\Omega_t)) := \left\{u \in L^p(0,T,L^q(\Omega_t)), \partial^\alpha u \in L^p(0,T,L^q(\Omega_t)) \ \ \forall \abs{\alpha} \leq l\right\}
\]
with the norm
\[
\norma{u}{L^p(0,T,W^{l,q}(\Omega_t))} := \sum_{\abs{\alpha} \leq l} \norma{\partial^\alpha u}{L^p(0,T,L^q(\Omega_t))}.
\]

The spaces of continuous functions in time with values in Lebesgue or Sobolev spaces in space variable $C([0,T];W^{l,q}(\Omega_t))$ are defined similarly as in \eqref{prostor} using the large ball $B_R(0)$.

We also recall that we denote as usual $H^k := W^{k,2}$. If no confusion may arise, we often drop the time interval and also the spatial domain to denote $L^p(W^{l,q}) := L^p(0,T;W^{l,q}(\Omega_t))$.

Finally we introduce a compact notation for the regularity class of the velocity in Theorem 1.1. For a function $f$ defined on $(0,T)\times\Omega_t$ we define
\begin{equation} \label{def:X}
\|f\|_{{\cal{X}}(T)}=\|f\|_{L^\infty(0,T;H^2(\Omega_t))\cap L^2(0,T;H^3(\Omega_t))} +\|f_t\|_{L^\infty(0,T;H^1(\Omega_t)) \cap L^2(0,T;H^2(\Omega_t))}+\|f_{tt}\|_{L^2(0,T;L^2(\Omega_t))}
\end{equation}
and for a function $\tilde f$ defined of $(0,T) \times \Omega$ where $\Omega$ is a fixed domain
\begin{equation} \label{def:Y}
\|\tilde f\|_{{\cal{Y}}(T)}=\|\tilde f\|_{L^\infty(0,T,H^2(\Omega))\cap L^2(0,T;H^3(\Omega))} +\|\tilde f_t\|_{L^\infty(0,T;H^1(\Omega)) \cap L^2(0,T;H^2(\Omega))}
+\|\tilde f_{tt}\|_{L^2(0,T;L^2(\Omega))}
\end{equation}
and obviously we denote by ${\cal X}(T)$ and ${\cal Y}(T)$ spaces for which above norms are finite.

\subsection{Weak solutions}
For weak solutions it is enough to assume the initial condition in a more general form
\begin{equation} \label{ic:weak}
\vr(0, \cdot) = \vr_0 ,\quad
(\vr\vu)(0, \cdot) = (\vr\vu)_0 \quad \mbox{in}\ \Omega_0.
\end{equation}
We define weak solutions to the compressible Navier-Stokes system with slip boundary conditions on moving domains as follows.
\begin{Definition}\label{d:ws}
We say that the couple $(\vr,\vu)$ is a weak solution of problem \eqref{i1a}-\eqref{i1b} with boundary conditions \eqref{i6b}-\eqref{i6c} and initial conditions \eqref{ic:weak} if
\begin{itemize}
\item $\vr \in L^\infty(0,T;L^\gamma(\mathbb{R}^3))$, $\vr \geq 0$ a.e. in $Q_T$ where $\gamma$ is from Theorem \ref{t:WSU}.
\item $\vu,\nabla_x\vu \in L^2(Q_T)$, $(\vu - \vc{V}) \cdot \vc{n}  (\tau , \cdot)|_{\Gamma_\tau}  = 0$ for a.a. $\tau \in [0,T]$.
\item The continuity equation \eqref{i1a} is satisfied in the whole $\mathbb{R}^3$ provided the density is extended by zero outside $\Omega_t$, i.e.
\bFormula{m1}
\int_{\Omega_\tau} \vr \varphi (\tau, \cdot) \ \dx - \int_{\Omega_0} \vr_0 \varphi (0, \cdot) \ \dx =
\int_0^\tau \int_{ \Omega_t} \left( \vr \partial_t \varphi + \vr \vu \cdot \Grad \varphi \right) \ \dxdt
\eF
for any $\tau \in [0,T]$ and any test function $\varphi \in \DC([0,T] \times \mathbb{R}^3)$.
\item Moreover, equation \eqref{i1a} is also satisfied in the sense of
renormalized solutions introduced by DiPerna and Lions \cite{DL}:
\bFormula{m2}
\int_{\Omega_\tau} b(\vr) \varphi (\tau, \cdot) \ \dx - \int_{\Omega_0} b(\vr_0) \varphi (0, \cdot) \ \dx =
\int_0^\tau \int_{ \Omega_t} \left( b(\vr) \partial_t \varphi + b(\vr) \vu \cdot \Grad \varphi +
\left( b(\vr)  - b'(\vr) \vr \right) \Div \vu \varphi \right) \ \dxdt
\eF
for any $\tau \in [0,T]$, any $\varphi \in \DC([0,T] \times \mathbb{R}^3)$, and any $b \in C^1 [0, \infty)$, $b(0) = 0$, $b'(r) = 0$ for large $r$.
\item The momentum equation \eqref{i1b} is replaced by the family of integral identities
\bFormula{m3}
\int_{\Omega_\tau} \vr \vu \cdot \vph (\tau, \cdot) \ \dx - \int_{\Omega_0} (\vr \vu)_0 \cdot \vph (0, \cdot) \ \dx
\eF
\[
= \int_0^\tau \int_{\Omega_t} \left( \vr \vu \cdot \partial_t \vph + \vr [\vu \otimes \vu] : \Grad \vph + p(\vr) \Div \vph
- \tn{S} (\Grad \vu) : \Grad \vph\right) \dxdt
-\kappa\int_0^\tau \int_{\Gamma_{\tau}}(\vu-\vV)\cdot\vph d\sigma
\]
for any $\tau \in [0,T]$ and any test function $\vph \in \DC([0,T] \times \mathbb{R}^3)$ satisfying
\bFormula{m4}
\vph \cdot \vn|_{\Gamma_\tau} = 0 \ \mbox{for any} \ \tau \in [0,T].
\eF
\end{itemize}
\end{Definition}

The existence of weak solutions to the problem \eqref{i1a}-\eqref{i7} was proved in \cite{FKNNS}.
The definition of weak solution in \cite{FKNNS} does not contain the boundary term since authors made assumptions  that $\kappa$ the friction coefficient  is zero. However, for the positive friction this term always appear with a good sign in the estimates. This means that it can even help for getting good a priori  estimates.  
\begin{Theorem}\label{t:ex}
Let $\Omega_0 \subset \mathbb{R}^3$ be a bounded domain of class $C^{2 + \nu}$. Assume $\vV$ and $p(\vr)$ satisfy the assumptions of Theorem \ref{t:WSU}.
Let the initial data \eqref{ic:weak} satisfy
\[
\vr_0 \in L^\gamma (\mathbb{R}^3),\ \vr_0 \geq 0, \ \vr_0 \not\equiv 0,\ \vr_0|_{\mathbb{R}^3 \setminus \Omega_0} = 0,\
(\vr \vu)_0 = 0 \ \mbox{a.a. on the set} \ \{ \vr_0 = 0 \} ,\ \int_{\Omega_0} \frac{1}{\vr_0} |(\vr \vu)_0 |^2 \ \dx < \infty.
\]

Then the problem \eqref{i1a}-\eqref{i6c},\eqref{ic:weak} admits a weak solution on any time interval $(0,T)$ in the sense specified through Definition \ref{d:ws}.
\end{Theorem}

We also notice that the existence theorem for the problem with no-slip boundary condition was proved in \cite{FeNeSt}.

\subsection{Iterative scheme}\label{ss:scheme}
Theorem \ref{Tmain} will be proved with the method of successive approximations.
At each step we will solve the linear system solving coupled linear continuity and
momentum equations. Here we adapt nowadays classical approach (\cite{V1}, \cite{Za}).
The linear continuity equation is solved in a moving domain which is possible since
the characteristics are well defined due to boundary condition \eqref{i6b}.
Then in order to solve the momentum equation we use Lagrangian transformation determined
by the velocity field $\vV$. It should be noticed that such approach is admissible since
the transformation depends only on $\vV$ and not on the solution, therefore is independent
on the step of iteration.
Then the crucial difficulty is in showing appropriate estimates that will give
convergence of the iterative scheme. We restrict our presentation to the proof
of the result for slip boundary conditions. With our method this case is more complicated
as it requires treatment of boundary terms which otherwise vanish.

In view of the above considerations we define our iterative scheme as follows.
We set $\rho_1(t,\vc{X}(t,x)) := \rho_0(x)$ and $\vc{u}_1(t,\vc{X}(t,x)) := \vc{u}_0(x)$ for all $t \in (0,T)$.

Assume we already have $(\vr_n,\vu_n)$. We define the next step of approximation $(\vr_{n+1},\vu_{n+1})$ as follows.


{\bf 1.} We solve for $\vr_{n+1}$ the linear continuity equation
\begin{equation} \label{eq:CEn}
\de_t\vr_{n+1}+\vu_n \cdot\nabla_x \vr_{n+1}
+\vr_{n+1}\div_x \vu_n=0 \quad \textrm{in}\; Q_T
\end{equation}
with the initial condition $\vr_{n+1}(0,x) = \vr_0(x)$ in $\Omega_0$.

{\bf 2.}
We solve for $\vu_{n+1}$ the linear momentum equation
\begin{align}\label{eq:MEn}
\vr_{n+1} \de_t \vu_{n+1}
- \mu\Delta_x\vu_{n+1} - (\frac{\mu}{3} + \eta)\nabla_x{\rm div}_x\vu_{n+1}= \\ \nonumber
- \vr_{n+1} \vu_{n} \cdot \nabla_x \vu_{n} - \nabla_x p(\vr_{n+1}) =: \vF(\vr_{n+1},\vu_n) \quad \textrm{in}\; Q_T
\end{align}
with boundary conditions
\begin{equation}
(\vu_{n+1}-\vc{V})\cdot\vc{n}|_{\Gamma_t}=0
\end{equation}
\begin{equation}
[\tn{S}(\nabla_x \vu_{n+1})\vn]_{\rm tan} + \kappa [\vu_{n+1}-\vV]_{\rm tan}|_{\Gamma_t}=0
\end{equation}
for all $t \in (0,T)$ and initial condition $\vu_{n+1}(0,x) = \vu_0(x)$ in $\Omega_0$.

\section{Linear continuity equation}\label{s:LinCont}

In order to have the iterative scheme well defined we have to solve in particular the linear continuity equation
\begin{equation} \label{lin_CE}
\vr_t+\vv \cdot \nabla_x \vr + \vr \div_x \vv=0  \quad \textrm{in}\; Q_T
\end{equation}
with $(\vv - \vV) \cdot \vc{n}|_{\Gamma_t}=0$. As pointed out earlier, we treat the linear continuity equation \eqref{lin_CE} directly in the moving domain and we do not use here any change of variables. The following result gives existence of solution to \eqref{lin_CE}.
\bProposition{sol_cont}
Assume $\vr_0 \in H^2(\Omega_0)$,
$\vv \in L^\infty(0,T;H^2(\Omega_t)) \cap L^2(0,T;H^3(\Omega_t))$.
Then for $T>0$ sufficiently small there exists a unique
solution $\vr$ to the linear transport equation \eqref{lin_CE} such that
\bFormula{sol_lin_cont}
\vr \in C([0,T];H^2(\Omega_t)), \quad \de_t \vr \in L^2(0,T;H^1(\Omega_t)).
\eF
Moreover the following estimates hold
\begin{equation} \label{est_lin_cont1}
\|\vr\|_{L^\infty(H^2)} \leq
C \|\vr_0\|_{H^2} \phi (\sqrt{T}\|\vv\|_{L^2(H^3)})=:D_0
\end{equation}
and
\begin{equation} \label{est_lin_cont2}
\|\de_t \vr\|_{L^2(H^1)} \leq C\sqrt{T} \|\vr_0\|_{H^2} \phi (\sqrt{T}\|\vv\|_{L^2(H^3)})\|\vv\|_{L^2(H^3)}
\end{equation}
where $\phi(\cdot)$ is an increasing, positive Lipschitz function.
Moreover, if $\vr_0\geq \kappa>0$ then $\vr \geq C(\kappa,T,\vv)>0$.
\eP
{\bf Remark.} In what follows we will denote by $\phi(\cdot)$ an increasing, positive,
sufficiently smooth function which may vary from line to line.\\

In order to obtain the desired form of estimates with the factor $\sqrt{T}$ we
use the solution formula provided by the method of characteristics. We have
\bFormula{formula_rho}
\vr(t,X(t,z))=\vr_0(z){\rm exp}\big(-\int_0^t \div_x{\vv}(s,X(s,z)){\rm d}s\big),
\eF
where
\begin{equation} \label{char}
X(t,z)=z+\int_0^t \vv(s,X(s,z)){\rm d}s.
\end{equation}
Note that the mapping $X(t,z)$ is different from the mapping $\vc{X}(t,z)$ introduced in \eqref{eq:coc}, since one is related to the velocity field $\vc{v}$ while the other to the velocity field $\vc{V}$. Before we proceed with the proof of Proposition \ref{Psol_cont} we first need some properties of the mapping $X$ and its derivatives.

\bLemma{dif_change}
Let $\vv \in L^2(0,T,H^3(\Omega_t))$. Let $X(t,z)$ be defined by \eqref{char}, i.e. for fixed $t \in (0,T)$ we have $X(t,\cdot): \Omega_0 \sil \Omega_t$. Then there exists sufficiently small $T > 0$ such that for all $t \in (0,T)$ there exists an inverse mapping $z(t,\cdot): \Omega_t \sil \Omega_0$, i.e. $z(t,X(t,y)) = y$ for all $y \in \Omega_0$ and $X(t,z(t,y)) = y$ for all $y \in \Omega_t$. Moreover it holds
\begin{equation} \label{nabla_x_1}
\|\nabla_x z(t,x) - \mathbb{I}\|_{L^\infty(Q_\tau)} \leq E(\tau)
\end{equation}
and
\begin{equation} \label{nabla_x_2}
\|\nabla^2_z X(t,z)\|_{L^\infty(0,T,L^4(\Omega_0))} \leq \phi(\|\vv\|_{L^2(0,T,H^3(\Omega_t))}),
\end{equation}
where $\mathbb{I}$ denotes the identity matrix, $E(t)$ is a nonnegative function such that $E(t) \sil 0$ as $t \sil 0^+$ and $\phi$ is an increasing positive function.
\eL

\bProof[Lemma \ref{Ldif_change}]
We have
\begin{equation} \label{dzx_1}
\nabla_z X(t,z)= \mathbb{I} +\int_0^t \nabla_x \vv(s,X(s,z)) \nabla_z X(s,z)ds,
\end{equation}
hence $\nabla_z X$ satisfies a system of ODE
$$
\de_t \nabla_z X = \nabla_x \vv(t,X(t,z))\nabla_z X.
$$
Multiplying the component of the above equation corresponding to $\de_{z_j}X_i$ by
$|\de_{z_j}X_i|^{p-2}\de_{z_j}X_i$ and integrating over $\Omega_0$ we obtain
$$
\frac{d}{d t}\|\nabla_z X\|_{L^p(\Omega_0)}^p \leq Cp \|\nabla_x \vv\|_{H^2(\Omega_t)} \|\nabla_z X\|_{L^p(\Omega_0)}^p,
$$
which can be simplified to 
$$
\frac{d}{d t}\|\nabla_z X\|_{L^p(\Omega_0)} \leq C \|\nabla_x \vv\|_{H^2(\Omega_t)} \|\nabla_z X\|_{L^p(\Omega_0)}.
$$
Therefore by Gronwall's inequality
$$
{\rm sup}_{0 < t < T} \|\nabla_z X\|_{L_p(\Omega_0)}\leq C {\rm exp}\left(\int_0^T\|\vv\|_{H^3(\Omega_s)}ds\right)
$$
which yields 
\begin{equation} \label{nabla_x_0}
\|\nabla_z X\|_{L^\infty(0,T,L^p(\Omega_0))} \leq C {\rm exp}\left(\int_0^T \|\vv\|_{H^3(\Omega_t)}ds\right)\leq \phi(\|\vv\|_{L^2(H^3)})\leq M.
\quad \forall \; 1 \leq p < \infty.
\end{equation}
Notice that the above estimate is uniform in $p$,
therefore we can justify the limit passage $p \sil \infty$ to conclude that \eqref{nabla_x_0} holds also for $p = \infty$.
In order to show the bound for the second
gradient we differentiate \eqref{dzx_1} w.r.t. $z$ and $t$ obtaining
$$
\frac{\de}{\de t} \nabla^2_z X(t,z) \sim \nabla^2_x \vv(t,X(t,z))(\nabla_z X(t,z))^2 + \nabla_x \vv(t,X(t,z))\nabla^2_z X(t,z).
$$
Note that we don't really care about the precise structure of the right hand side, in order to obtain estimates, it is enough for us to know that we have term involving second gradient of $\vc{v}$ multiplied twice by first gradient of $X$ and second gradient of $X$ multiplied by first gradient of $\vc{v}$.

We proceed similarly as before. Multiplying
the component corresponding to $\nabla^2_{z_i z_j}X_k$ by $|\nabla^2_{z_i z_j}X_k|^2\nabla^2_{z_i z_j}X_k$
and integrating over $\Omega_0$ we obtain
$$
\frac{d}{d t}\|\nabla^2_z X(t,\cdot)\|_{L^4}^4 \leq C
\int_{\Omega_0} |\nabla_x^2 \vv(t,X(t,\cdot))| |\nabla_z X(t,\cdot)|^2 |\nabla_z^2 X(t,\cdot)|^3\dx + C
\int_{\Omega_0} |\nabla_x \vv(t,X(t,\cdot))| |\nabla_z^2 X(t,\cdot)|^4\dx \leq
$$$$
\leq C\|\nabla^2_x \vv(t)\|_{L^6(\Omega_t)}\|\nabla_z X(t)\|^2_{L^{24}(\Omega_0)}\|\nabla_z^2 X(t)\|_{L^4(\Omega_0)}^3
+\|\nabla_x \vv(t)\|_{H^2(\Omega_t)}\|\nabla_z^2 X(t)\|_{L^4(\Omega_0)}^4,
$$
which by Sobolev embedding and Gronwall inequality implies \eqref{nabla_x_2}.

Now we are ready to finish the proof of \eqref{nabla_x_1}. Recall that from \eqref{nabla_x_0} we have
\begin{equation} \label{dzx_2}
\|\nabla_z X(t,z)\|_{L^\infty((0,T) \times \Omega_0)} \leq M
\end{equation}
for some $M > 0$ which depends on the $L^2(H^3)$ norm of $\vv$. Therefore from \eqref{dzx_1} we have
\begin{equation} \label{nabla_x_3}
|\nabla_z X(t,z)-\mathbb{I}|=\left|\int_0^t \nabla_x \vv(s,X(s,z))\nabla_z X(s,z)ds\right| \leq M\norma{\nabla_x\vv}{L^1(0,\tau;L^\infty(\Omega_t))} \leq CM\sqrt{\tau}\norma{\vv}{L^2(0,\tau;H^3(\Omega_t))}.
\end{equation}
The expression on the right hand side of \eqref{nabla_x_3} will be small for small times, which yields the
invertibility of $\nabla_z X$ and also the bound \eqref{nabla_x_1}.

\qed

Now we proceed with the proof of Proposition \ref{Psol_cont}.

\bProof[Proposition \ref{Psol_cont}]
The solution formula \eqref{formula_rho} immediately gives the last statement of the lemma.
Moreover, denoting $\bar \vr(t,z) = \vr(t,X(t,z))$ we obtain from \eqref{formula_rho}
\begin{equation} \label{rho00}
\| \vr(t,\cdot)\|_{L^2(\Omega_t)} \leq C\|\bar \vr(t,\cdot)\|_{L^2(\Omega_0)} \leq C\|\vr_0\|_{L^2} \|{\rm exp}\big(-\int_0^t \ddiv_x \vv \d s\big)\|_{L^\infty}\leq C\|\vr_0\|_{L^2} {\rm exp}\big(\int_0^t \|\ddiv_x \vv \|_{L^\infty} \d s\big)
\end{equation}
and consequently
\begin{equation} \label{rho001}
\| \vr\|_{L^\infty(0,T;L^2(\Omega_t))} \leq C\|\bar \vr\|_{L^\infty(0,T;L^2(\Omega_0))} \leq C\|\vr_0\|_{L^2} {\rm exp}\big(\|\ddiv_x \vv \|_{L^1(0,T,L^\infty)} \big) \leq \|\vr_0\|_{L^2}\phi(\sqrt{T}\|\vv \|_{L^2(0,T;H^3)}).
\end{equation}

%
Differentiating \eqref{formula_rho} w.r.t $z$ we obtain
\begin{align}\label{nablarho}
&\nabla_x\vr(t,X(t,z))\nabla_z X(t,z)  \\ \nonumber
&\qquad = {\rm exp}(-\int_0^t \div_x \vv(s,X(s,z))\, {\rm d} s)\left(\nabla_z\vr_0(z) - \vr_0(z)\int_0^t\nabla_x\div_x\vv(s,X(s,z))\nabla_z X(s,z)\,{\rm d}s\right)
\end{align}
and thus at least for small times
\begin{equation} \label{grad_vr}
|\nabla_x \vr| \leq |(\nabla_z X)^{-1}| \left( \big| \nabla_z \vr_0\, {\rm exp}(-\int_0^t \ddiv_x \vv)\big|
+ \big|\vr_0\, {\rm exp}(-\int_0^t \ddiv_x \vv) \int_0^t \nabla_x \ddiv_x \vv \nabla_z X\, {\rm d}s\big| \right).
\end{equation}
Now, the first term in \eqref{grad_vr} can be estimated as before (see \eqref{rho00} and \eqref{rho001}) and with the second term we proceed as follows
\begin{align} \label{grad_vr2}
&\|(\nabla_z X)^{-1}\, \vr_0\, {\rm exp}(-\int_0^t \ddiv_x \vv) \int_0^t \nabla_x \ddiv_x \vv \nabla_z X\,{\rm d}s\big\|_{L^2(\Omega_0)} \leq \\ \nonumber
&\|(\nabla_z X)^{-1}\|_{L^\infty(\Omega_0)}\|\vr_0\|_{L^\infty(\Omega_0)}\|{\rm exp}(-\int_0^t \ddiv_x \vv)\|_{L^\infty(\Omega_0)} \int_0^t (\|\nabla_x \ddiv_x \vv\|_{L^2(\Omega_t)}\|\nabla_z X\|_{L^\infty(\Omega_0)})\,{\rm d}s \leq \\ \nonumber
&\|(\nabla_z X)^{-1}\|_{L^\infty(\Omega_0)}\|\vr_0\|_{H^2(\Omega_0)}{\rm exp}(\int_0^t \|\ddiv_x \vv\|_{L^\infty(\Omega_t)}) \int_0^t (\|\nabla_x \ddiv_x \vv\|_{L^2(\Omega_t)}\|\nabla_z X\|_{L^\infty(\Omega_0)})\,{\rm d}s.
\end{align}
Taking supremum over $t \in (0,T)$, using the fact that both $\|(\nabla_z X)^{-1}\|_{L^\infty(\Omega_0)}$ and $\|\nabla_z X\|_{L^\infty(\Omega_0)}$ are bounded in time for sufficiently small $T$ and using H\"older inequality similarly as in \eqref{rho001} to obtain the factor $\sqrt{T}$ we arrive at
\begin{equation} \label{est_grad_vr}
\|\nabla_x \vr\|_{L^\infty(0,T;L^2(\Omega_t))}\leq C\|\vr_0\|_{H^2(\Omega_0)}\phi(\sqrt{T}\|\vv\|_{L^2(0,T;H^3(\Omega_t))}).
\end{equation}

In order to estimate the second derivatives of $\vr$ we differentiate \eqref{nablarho} one more time with respect to $z$. Again we are not particularly interested in the precise structure of the resulting equation, for the purpose of the estimate it is enough to recognize all kinds of terms appearing on both sides of the resulting equation. We have
\begin{align}\label{nabla2rho}
&\nabla^2_x \vr (\nabla_z X)^2 + \nabla_x \vr \nabla_z^2 X  \sim \\ \nonumber
& \qquad \sim \nabla^2_z\vr_0\,{\rm exp}(-\int_0^t \ddiv_x \vv) - 2\nabla_z\vr_0\,{\rm exp}(-\int_0^t \ddiv_x \vv)\int_0^t(\nabla_x\div_x\vv \nabla_z X)\, {\rm d}s \\ \nonumber
& \qquad + \vr_0\,{\rm exp}(-\int_0^t \ddiv_x \vv)\left(\int_0^t(\nabla_x\div_x\vv \nabla_z X)\, {\rm d}s\right)^2 - \vr_0\,{\rm exp}(-\int_0^t \ddiv_x \vv)\int_0^t(\nabla_x\div_x\vv \nabla^2_z X)\, {\rm d}s \\ \nonumber
& \qquad - \vr_0\,{\rm exp}(-\int_0^t \ddiv_x \vv)\int_0^t(\nabla^2_x\div_x\vv (\nabla_z X)^2)\, {\rm d}s,
\end{align}
therefore 
\begin{equation}
\nabla^2_x \vr \sim [(\nabla_z X)^2]^{-1}\left[ -\nabla_x \vr \nabla_z^2 X + \textrm{RHS of \eqref{nabla2rho}} \right].
\end{equation}
Since by \eqref{nabla_x_1} $[(\nabla_z X)^2]^{-1}$ is bounded for small times, it remains to estimate 
the $L^\infty(L^2)$ norm of the right hand side of \eqref{nabla2rho} and also of $\nabla_x \vr \nabla_z^2 X$. Let us start with the latter. Due to Lemma \ref{Ldif_change} we have the $L^\infty(L^4)$ bound for $\nabla_z^2 X$ and for $\nabla_x \vr$ we use the formula \eqref{grad_vr}. We show only the estimates in the space variable, the time variable is handled always the same way as in \eqref{rho001} or \eqref{est_grad_vr}. We have
\begin{align}
\norma{\nabla_z^2 X(\nabla_z X)^{-1}{\rm exp}(-\int_0^t \ddiv_x \vv)\nabla_z \vr_0}{L^2} &\leq \norma{\nabla_z^2 X}{L^4}\norma{(\nabla_z X)^{-1}}{L^\infty}\norma{{\rm exp}(-\int_0^t \ddiv_x \vv)}{L^\infty}\norma{\nabla_z \vr_0}{L^4} \\ \nonumber
& \leq C\norma{\nabla_z^2 X}{L^4}\norma{(\nabla_z X)^{-1}}{L^\infty}{\rm exp}(\int_0^t\norma{\ddiv_x \vv}{L^\infty} )\norma{\vr_0}{H^2}
\end{align}
and
\begin{align}\label{similar}
&\norma{\nabla_z^2 X(\nabla_z X)^{-1}{\rm exp}(-\int_0^t \ddiv_x \vv)\vr_0\int_0^t(\nabla_x\ddiv_x \vv \nabla_z X)}{L^2} \\ \nonumber
&\qquad \leq \norma{\nabla_z^2 X}{L^4}\norma{(\nabla_z X)^{-1}}{L^\infty}\norma{{\rm exp}(-\int_0^t \ddiv_x \vv)}{L^\infty}\norma{\vr_0}{L^\infty}\int_0^t(\norma{\nabla^2_x \vv}{L^4}\norma{\nabla_z X}{L^\infty}) \\ \nonumber
& \qquad \leq C\norma{\nabla_z^2 X}{L^4}\norma{(\nabla_z X)^{-1}}{L^\infty}{\rm exp}(\int_0^t\norma{\ddiv_x \vv}{L^\infty} )\norma{\vr_0}{H^2}\int_0^t(\norma{\vv}{H^3}\norma{\nabla_z X}{L^\infty}).
\end{align}

Now let us treat the right hand side of \eqref{nabla2rho}. The first term can be handled similarly as above in \eqref{rho001}. For the second term we have
\begin{align}
\norma{\nabla_z\vr_0{\rm exp}(-\int_0^t \ddiv_x \vv)\int_0^t(\nabla_x\ddiv_x \vv \nabla_z X)}{L^2} & \leq \norma{\nabla_z \vr_0}{L^4}\norma{{\rm exp}(-\int_0^t \ddiv_x \vv)}{L^\infty}\int_0^t(\norma{\nabla^2_x \vv}{L^4}\norma{\nabla_z X}{L^\infty}) \\ \nonumber
& \leq C\norma{\vr_0}{H^2}{\rm exp}(\int_0^t\norma{\ddiv_x \vv}{L^\infty} )\int_0^t(\norma{\vv}{H^3}\norma{\nabla_z X}{L^\infty}).
\end{align}
The third term on the right hand side of \eqref{nabla2rho} is estimated as
\begin{align}
&\norma{\vr_0\,{\rm exp}(-\int_0^t \ddiv_x \vv)\left(\int_0^t(\nabla_x\div_x\vv \nabla_z X)\, {\rm d}s\right)^2}{L^2} \\ \nonumber
& \qquad \leq \norma{\vr_0}{L^\infty}\norma{{\rm exp}(-\int_0^t \ddiv_x \vv)}{L^\infty}\left(\int_0^t(\norma{\nabla^2_x \vv}{L^4}\norma{\nabla_z X}{L^\infty})\right)^2 \\ \nonumber
& \qquad \leq C\norma{\vr_0}{H^2}{\rm exp}(\int_0^t\norma{\ddiv_x \vv}{L^\infty} )\left(\int_0^t(\norma{\vv}{H^3}\norma{\nabla_z X}{L^\infty})\right)^2.
\end{align}
The fourth term is estimated similarly as in \eqref{similar}. We have
\begin{align}
&\norma{\vr_0\,{\rm exp}(-\int_0^t \ddiv_x \vv)\int_0^t(\nabla_x\div_x\vv \nabla^2_z X)}{L^2} \\ \nonumber
&\qquad \leq \norma{\vr_0}{L^\infty}\norma{{\rm exp}(-\int_0^t \ddiv_x \vv)}{L^\infty}\int_0^t(\norma{\nabla^2_x \vv}{L^4}\norma{\nabla_z^2 X}{L^4}) \\ \nonumber
& \qquad \leq C\norma{\vr_0}{H^2}{\rm exp}(\int_0^t\norma{\ddiv_x \vv}{L^\infty} )\int_0^t(\norma{\vv}{H^3}\norma{\nabla_z^2 X}{L^4}).
\end{align}
The last term on the right hand side of \eqref{nabla2rho} is the only one involving third derivatives of the velocity field $\vv$. We proceed as follows
\begin{align}
&\norma{\vr_0\,{\rm exp}(-\int_0^t \ddiv_x \vv)\int_0^t(\nabla^2_x\div_x\vv (\nabla_z X)^2)}{L^2} \\ \nonumber
&\qquad \leq \norma{\vr_0}{L^\infty}\norma{{\rm exp}(-\int_0^t \ddiv_x \vv)}{L^\infty}\int_0^t(\norma{\nabla^3_x \vv}{L^2}\norma{\nabla_z X}{L^\infty}^2) \\ \nonumber
& \qquad \leq C\norma{\vr_0}{H^2}{\rm exp}(\int_0^t\norma{\ddiv_x \vv}{L^\infty} )\int_0^t(\norma{\vv}{H^3}\norma{\nabla_z X}{L^\infty}^2)
\end{align}
As we mentioned above, handling the time variable is the same as in the estimates for $\vr$ and $\nabla_x \vr$. Since we have the $L^2(H^3)$ bound for $\vv$ and under the time integrals there is always only first power of the $H^3$ norm of $\vv$ we always gain the factor $\sqrt{T}$ and finally end up with
\begin{equation} \label{est_grad2_vr}
\|\nabla_x^2 \vr\|_{L^\infty(0,T;L^2(\Omega_t))}\leq C\|\vr_0\|_{H^2(\Omega_0)}\phi(\sqrt{T}\|\vv\|_{L^2(0,T;H^3(\Omega_t))}).
\end{equation}
Combining \eqref{rho00}, \eqref{est_grad_vr} and \eqref{est_grad2_vr} we obtain \eqref{est_lin_cont1}. 
The continuity in time of the $H^2$ norm of $\vr$ follows from the solution formula \eqref{formula_rho}.

%

With the estimates of $\vr$ in $L^\infty(H^2)$ at hand we can use the equation \eqref{lin_CE} to obtain estimates for the time derivative of $\rho$. We have
\begin{equation}\label{eq:bla}
\de_t \vr = - \vv \cdot \nabla_x \vr - \vr \ddiv_x \vv.
\end{equation}
Estimating the right hand side as
\begin{align}
\norma{\de_t \vr}{L^2(0,T;L^2(\Omega_t))} &\leq \norma{\vc{v}\cdot\nabla_x \vr}{L^2(0,T;L^2(\Omega_t))} + \norma{\vr\div_x\vc{v}}{L^2(0,T;L^2(\Omega_t))} \\ \nonumber
&\leq C\sqrt{T}\left(\norma{\vc{v}}{L^\infty(0,T;H^2(\Omega_t))}\norma{\nabla_x\vr}{L^\infty(0,T;L^2(\Omega_t))} + \norma{\vr}{L^\infty(0,T;H^2(\Omega_t))}\norma{\div_x \vc{v}}{L^\infty(0,T;L^2(\Omega_t))}\right) \\ \nonumber
&\leq C\sqrt{T}\|\vr_0\|_{H^2(\Omega_0)}\phi(\sqrt{T}\|\vv\|_{L^2(0,T;H^3(\Omega_t))})\norma{\vc{v}}{L^\infty(0,T;H^2(\Omega_t))}.
\end{align}

Next we apply the spatial gradient to \eqref{eq:bla} to obtain
\begin{equation}
\de_t \nabla_x \vr = - \vv \cdot \nabla_x \nabla_x\vr - \nabla_x \vr\nabla_x \vv - \nabla_x \vr \ddiv_x \vv - \vr\nabla_x\ddiv_x \vv
\end{equation}
and we estimate the $L^2(L^2)$ norm of the right hand side. Starting first with the norms in the $x$ variable we have
\begin{align}
\norma{\vc{v} \nabla^2_x \vr}{L^2(\Omega_t)} &\leq \norma{\vc{v}}{L^\infty}\norma{\nabla^2_x\vr}{L^2} \leq C \norma{\vc{v}}{H^2}\norma{\vr}{H^2} \\
\norma{\nabla_x \vr \nabla_x \vc{v}}{L^2(\Omega_t)} &\leq \norma{\nabla_x\vr}{L^4}\norma{\nabla_x\vc{v}}{L^4} \leq C \norma{\vc{v}}{H^2}\norma{\vr}{H^2} \\
\norma{\vr \nabla^2_x \vc{v}}{L^2(\Omega_t)} &\leq \norma{\vr}{L^\infty}\norma{\nabla^2_x\vc{v}}{L^2} \leq C \norma{\vc{v}}{H^2}\norma{\vr}{H^2},
\end{align}
so altogether we conclude
\begin{align}
\norma{\de_t \nabla_x\vr}{L^2(0,T,L^2(\Omega_t))} &\leq C\sqrt{T}\norma{\vc{v}}{L^\infty(H^2)}\norma{\vr}{L^\infty(H^2)} \\ \nonumber
&\leq C\sqrt{T}\|\vr_0\|_{H^2(\Omega_0)}\phi(\sqrt{T}\|\vv\|_{L^2(0,T;H^3(\Omega_t))})\norma{\vc{v}}{L^\infty(0,T;H^2(\Omega_t))}.
\end{align}

\qed

\section{Linear momentum equation}\label{s:LinMom}

In this section we treat the linear momentum equation
\begin{align} \label{lin_ME}
\vr \de_t \vu
- \mu\Delta_x\vu - (\frac{\mu}{3} + \eta)\nabla_x\div_x\vu = \vc{F} \quad \textrm{in}\; Q_T,\\
(\vu-\vc{V})\cdot\vc{n}|_{\Gamma_\tau}=0,\nonumber\\
\vn \cdot \mu(\nabla_x \vu + \nabla_x^t \vu) \cdot \tau_k + \kappa (\vu-\vV)\cdot \tau_k|_{\Gamma_\tau}=0, \quad k=1,2.
\end{align}

The next proposition gives existence of solutions to the linear momentum equation on moving domain and we deliberately skip emphasizing the domains of all function spaces in this proposition in order to shorten the notation. However we recall here that by $L^p(B)$ we mean in this lemma the function space $L^p(0,T;B(\Omega_t))$. We recall that the space ${\cal X}(T)$ was defined in \eqref{def:X}.
\bProposition{sol_mom}
Let $T > 0$ be sufficiently small. Let $\vV,p(\vr)$ satisfy the assumptions of Theorem \ref{Tmain}. Assume $\vr \in L^\infty(H^2)$, $\vr_t \in L^2(H^1)$,
$\vF \in L^2(H^1)$, $\vF_t \in L^2(L^2)$ and $\vF(0)\in H^1(\Omega_0)$, $\vu_0 \in H^3(\Omega_0)$.
Then there exist a unique solution $\vu$ to the system \eqref{lin_ME} such that
$\vu \in \cal{X}(T)$
and the following estimate holds
\begin{align} \label{est_lin_ME}
\|\vu\|_{{\cal{X}}(T)} &\leq \phi(\|\vr\|_{L^\infty(H^2)},\|\vr_t\|_{L^2(H^1)})  \\ \nonumber
&\times \left(\|\vF\|_{L^2(H^1)}+ \|\vF_t\|_{L^2(L^2)}+\|\vF(0)\|_{H^1}+\|\vu_0\|_{H^3}+\|\vV\|_{L^\infty(W^{1,\infty})}+\|\vV_t\|_{L^2(H^1)}+\|\vV_{tt}\|_{L^2(L^2)})\right),
\end{align}
where $\phi$ is a positive increasing function of its arguments.
\eP
This time it is more convenient to convert the problem to fixed spatial domain.
For this purpose we introduce the Lagrangian coordinates determined by $\vV$.
As this is an important and independent step we present it in a separate subsection.

\subsection{Lagrangian coordinates and linearization} \label{sec_lin}
We start with rewriting the problem \eqref{lin_ME} defined on $Q_T$ to a problem defined on a fixed spatial domain $(0,T) \times \Omega_0$ using time dependent change of coordinates. To this end we use the formula \eqref{eq:coc}. We set
\begin{equation}\label{eq:new}
\vrt(t,y) := \vr(t,\vc{X}(t,y)), \quad \vut(t,y) := \vu(t,\vc{X}(t,y))
\end{equation}
and we denote the components of the vector $\vc{X}$ as $\vc{X} = (X_1,X_2,X_3)$. To proceed we also need the inverse mapping to $\vc{X}(t,y)$ which we denote by $\vc{Y}(t,x)$, thus it holds for all $t \geq 0$ and all $x \in \Omega_t$
\begin{equation}\label{eq:inv}
\vc{X}(t,\vc{Y}(t,x)) = x
\end{equation}
and again we denote the components of $\vc{Y}$ as $\vc{Y} = (Y_1,Y_2,Y_3)$.

Differentiating \eqref{eq:inv} with respect to time we obtain
\begin{equation}\label{eq:td}
\frac{\de \vc{X}}{\de t} + \frac{\de \vc{X}}{\de y}\frac{\de \vc{Y}}{\de t} = 0
\end{equation}
and thus
\begin{equation}\label{eq:td2}
\frac{\de \vc{Y}}{\de t} = - \vV\cdot \Grad\vc{Y}.
\end{equation}
Using \eqref{eq:td2} we thus transform the time derivative of $u_i$ as follows
\begin{equation}\label{eq:tdr}
\frac{\de u_i}{\de t} = \frac{\de \ut_i}{\de t} + \nabla_y\ut_i \cdot \frac{\de \vc{Y}}{\de t} = \frac{\de \ut_i}{\de t} - \nabla_y\ut_i \cdot \left( \vV\cdot\nabla_x\vc{Y} \right).
\end{equation}
The spatial derivatives transform just by multiplying by the Jacobian of the change of coordinates.
Therefore the $i$-th component of the equation \eqref{lin_ME}$_1$
rewritten in terms of $\vrt,\vut$ defined on fixed domain $\Omega_0$ reads
%
\begin{align}\label{eq:mer1}
&\vrt\big[ \frac{\de \ut_i}{\de t} - \frac{\de \ut_i}{\de y_j} V_k \frac{\de Y_j}{\de x_k}\big]
- \mu \frac{\de^2 \ut_i}{\de y_k\de y_l}\frac{\de Y_l}{\de x_p}\frac{\de Y_k}{\de x_p}
\\ \nonumber
 - &\mu \frac{\de \ut_i}{\de y_k}\Delta_x Y_k - (\frac{\mu}{3}+\eta)\frac{\de^2 \ut_p}{\de y_k\de y_l}\frac{\de Y_l}{\de x_p}\frac{\de Y_k}{\de x_i} - (\frac{\mu}{3}+\eta)\frac{\de \ut_p}{\de y_k}\frac{\de^2 Y_k}{\de x_i\de x_p} = \tilde F_i,
\end{align}
where we have used Einstein summation convention.
This can be rewritten as a usual linearized momentum equation with a transport term on the left hand side and a right hand side containing terms which are either small for small times or of lower order:
\begin{equation}\label{eq:mer2}
\vrt\frac{\de \vut}{\de t}
- \mu\Delta_y(\vut) - (\frac{\mu}{3} + \eta)\nabla_y{\rm div}_y\vut
= \vF + \vrt\vV\cdot\nabla_y\vut + {\bf R}(\vrt,\vut)=:\tilde\vF(\vrt,\vut)
\end{equation}
with
\begin{align} \label{def_R}
{\bf R}(\vrt,\vut)=&\vrt\frac{\de \vut}{\de y_j} V_k\left(\frac{\de Y_j}{\de x_k} - \delta_{jk}\right)
+ \mu \frac{\de^2 \vut}{\de y_k\de y_l}\left(\frac{\de Y_l}{\de x_p}\frac{\de Y_k}{\de x_p} - \delta_{lp}\delta_{kp}\right) \\ \nonumber
&+ (\frac{\mu}{3}+\eta)\frac{\de^2 \ut_p}{\de y_k\de y_l}\left(\frac{\de Y_l}{\de x_p}\nabla_x Y_k - \delta_{lp}\vc{e}_k\right)
+ \mu \frac{\de \vut}{\de y_k}\Delta_x Y_k  + (\frac{\mu}{3}+\eta)\frac{\de \ut_p}{\de y_k}\nabla_x \frac{\de Y_k}{\de x_p},
\end{align}
where $\vc{e}_j$ is the $j$-th unit vector.

The boundary conditions are transformed in a nontrivial way as well. In particular, the condition \eqref{lin_ME}$_2$ transforms as
\begin{equation} \label{bc1}
(\vut - \vc{V})(t,y) \cdot \vc{n}(y) = (\vut - \vc{V})(t,y)\cdot\left(\vc{n}(y)-\vc{n}(\vc{X}(t,y))\right) + (\vc{V}(t,\vc{X}(t,y))-\vc{V}(t,y))\cdot \vc{n}(\vc{X}(t,y))=:\vc{d}(\vut,\vV)(t,y)
\end{equation}
for $y \in \partial\Omega_0$. Note however, that for small times, the expression $\vc{d}(\vut,\vV)$ will be small due to the fact that the mapping $\vc{X}$ is close to identity and its regularity is given by the regularity of $\vV$.

Similarly the boundary condition \eqref{lin_ME}$_3$ will also be transformed, however since it contains differentiation, the resulting expression is more complicated. We denote by $\tau^1,\tau^2$ tangent vectors to $\partial\Omega_0$. A lengthy yet straightforward computation yields for $y \in \partial\Omega_0$ the following:
\begin{align} \label{bc2}
&\mu(\nabla_y \vut+\nabla^T_y \vut)(t,y)\vc{n}(y) \cdot \tau^k(y) + \kappa (\vut-\vV)(t,y)\cdot\tau^k(y) = \\ \nonumber
& \qquad = \mu\left(\nabla_y\vut(t,y)(\tn{I}-\nabla_x\vc{Y}) + ((\tn{I}-\nabla_x^T\vc{Y})\nabla_y^T\vut(t,y))^T\right)\vc{n}(\vc{X}(t,y))\cdot\tau^k(\vc{X}(t,y)) \\ \nonumber
& \qquad + \mu(\nabla_y\vut+\nabla^T_y\vut)(t,y)[(\vc{n}(y)-\vc{n}(\vc{X}(t,y)))\cdot\tau^k(\vc{X}(t,y)) + \vc{n}(t,y)\cdot(\tau^k(y)-\tau^k(\vc{X}(t,y)))] \\ \nonumber
& \qquad + \kappa(\vut-\vV)(t,y)\cdot(\tau^k(y)-\tau^k(\vc{X}(t,y))) + \kappa (\vc{V}(t,\vc{X}(t,y))-\vc{V}(t,y))\cdot \tau^k(\vc{X}(t,y)) =: \vc{B}(\vut,\vV)(t,y).
\end{align}
Again we emphasize that despite the rather complicated structure of $\vc{B}(\vut,\vV)$ it is easy to observe, that this expression will be small for small times as a consequence of a fact that $\vc{X}$ is close to identity for small times.


The right hand side of \eqref{eq:mer2} and boundary conditions \eqref{bc1}, \eqref{bc2} contains the solution and variable coefficients dependent of the change of variables.
However, all these quantities will remain small for small times in appropriate norms.
Therefore what is important is the structure
of the left hand side and in particular we will be able to solve the system \eqref{eq:mer2}, \eqref{bc1}, \eqref{bc2} once we have solved the following linear
problem on a fixed domain $(0,T)\times\Omega_0$.
\begin{align} \label{lin_ME_fix}
\vr\frac{\de \vu}{\de t} 
- \mu\Delta_y(\vu) - (\frac{\mu}{3} + \eta)\nabla_y{\rm div}_y\vu
= \vF,\\
(\vu-\vV)\cdot\vc{n}|_{\Gamma_0}=\vc{d},\nonumber\\
\left(\left[ \tn{S}(\nabla_y\vu) \vc{n} \right]_{\rm tan} + \kappa\left[ \vu-\vV \right]_{\rm tan}\right)|_{\Gamma_0} = \vc{B}.\nonumber
\end{align}
For simplicity we denote the unknown of this system as $\vu$ instead of $\tilde{\vu}$. This system is solved in the next subsection.

\subsection{Solution of the linear momentum equation on a fixed domain}

In this subsection we work with a system of differential equations stated on a domain $(0,T) \times \Omega_0$. We again skip emphasizing the domains of all function spaces in this subsection in order to shorten the notation and we note here that by $L^p(B)$ we mean here the function space $L^p(0,T;B(\Omega_0))$. We also skip the index $y$ in differential operators.

\bLemma{sol_mom1} 
Assume that $\vc{B}$ and $\vc{d}$ admits an extension to $\Omega_0$ given by
\begin{align} \label{prob_ext}
\vu^b\cdot\vc{n}|_{\Gamma_0}=\vc{d},\nonumber\\
\left(\left[ \tn{S}(\nabla\vu^b) \vc{n} \right]_{\rm tan} + \kappa\left[ \vu^b \right]_{\rm tan}\right)|_{\Gamma_0} = \vc{B}
\end{align}
such that $\vu^b \in {\cal{Y}}(t)$, where ${\cal{Y}}(t)$ is defined in  \eqref{def:Y}. Let $p,\vV$ satisfy the assumptions of Theorem \ref{Tmain}.
Assume further that $\vr \in L^\infty(H^2)$, $\vr_t \in L^2(H^1)$, $\vr \geq C >0$,
$\vF \in L^2(H^1)$, $\vF_t \in L^2(L^2)$, $\vu_0 \in H^3(\Omega_0)$.
Then there exists a unique solution $\vu$ to the problem \eqref{lin_ME_fix}
such that
\begin{align} \label{est_lin_ME_fix}
&\|\vu\|_{{\cal{Y}}(T)} \leq \phi(\|\vr\|_{L^\infty(H^2)},\|\vr_t\|_{L^2(H^1)})
\left(\|\vF\|_{L^2(H^1)\cap L^\infty(L_2)}+ \|\vF_t\|_{L^2(L^2)} + \|\vu^b\|_{{\cal Y}(T)} \right. \\ \nonumber
& \qquad + \left.\|\vu_0\|_{H^3}+\|\vV\|_{L^\infty(W^{1,\infty})}+\|\vV_t\|_{L^2(H^1)}+\|\vV_{tt}\|_{L^2(L^2)})\right),
\end{align}
where
$\phi$ denotes a positive increasing function of its arguments.
\eL
\bProof
The result is based on Lemmas 2.2-2.4 from \cite{Za}
where the same linear system is considered, however with $\vc{d}=\vc{B}=\vV=0$.
Therefore we present here some details to show how we treat
the inhomogeneous boundary data. For simplicity of notation we set for the rest of this subsection $\Omega := \Omega_0$.

We start with removing the inhomogeneity from the boundary data.
For this purpose we take the extension $\vu^b$ defined by \eqref{prob_ext}.
Taking $\tilde \vu = \vu - \vu^b$ we obtain 
\begin{align} \label{lin_ME1}
\vr \de_t \vut
- \mu\Delta(\vut) - (\frac{\mu}{3} + \eta)\nabla\div\vut = \vc{\bar F}-\vr \de_t \vu^b,\\
(\vut-\vc{V})\cdot\vc{n}|_{\Gamma}=0,\nonumber\\
\left(\left[ \tn{S}(\nabla\vut) \vc{n} \right]_{\rm tan} + \kappa\left[ \vut - \vc{V} \right]_{\rm tan}\right)|_{\Gamma} = 0,\nonumber
\end{align}
where
\begin{equation} \label{def_bF}
\vc{\bar F} = \vF + \mu\Delta(\vu^b) + (\frac{\mu}{3} + \eta)\nabla\div\vu^b.
\end{equation}
We have to solve the above system, for simplicity we denote $\vut$ again by $\vu$.
We can also assume the friction coefficient
$\kappa=0$, positive friction yields only additional lower order terms which are easy to treat.

First we can write the weak formulation of \eqref{lin_ME1} and using $\vu-\vc{V}$ as a test function
(we have to test by a function with vanishing normal component at the boundary)
we obtain the energy estimate
\begin{equation} \label{ene}
\|\vu\|_{L^\infty(L^2)}+\|\nabla \vu\|_{L^2(L^2)} \leq C \|\vc{\bar F},\vV,\vu^b_t\|_{L^2(L^2)}.
\end{equation}
Next we multiply $\eqref{lin_ME1}$ by $\de_t (\vu-\vc{V}) + \ep A\vu$, where $\ep$
is sufficiently small and $A\vu=-\mu \Delta \vu - (\frac{\mu}{3}+\eta)\nabla \ddiv \vu$.
Integrating over $\Omega$ we get
\begin{equation} \label{est_lin3}
\int_{\Omega}\vr|\vu_t|^2\dx +\int_{\Omega}(\vu-\vc{V})_t \cdot A\vu\dx + \ep\int_{\Omega}|A\vu|^2\dx=
- \ep \int_{\Omega} \vr \vu_t\cdot A\vu\dx
+\int_{\Omega}\vr \vu_t\cdot \vc{V}_t \dx
+\int_{\Omega} (\vc{\bar F}-\vr \de_t \vu^b)\cdot(\vu_t+\ep A\vu)\dx
\end{equation}
Integrating by parts the second term on the left hand side we obtain
\begin{align}
\int_{\Omega}(\vu-\vc{V})_t\cdot A\vu\dx&=-\int_{\Omega}(\vu-\vc{V})_t\cdot\ddiv \tn{S}(\nabla \vu)\dx=
\int_{\Omega}\tn{S}(\nabla \vu):\nabla (\vu_t-\vV_t)\dx - \int_{\de \Omega}(\vu-\vc{V})_t \cdot \tn{S}(\nabla \vu)\vn\dS=\nonumber\\
&=\int_{\Omega}\tn{S}(\nabla \vu):\nabla \vu_t\dx-\int_{\Omega}\tn{S}(\nabla \vu):\nabla \vV_t\dx, \nonumber
\end{align}
where in the last step we decomposed $\tn{S}$ into tangential and normal part and used the boundary conditions. Notice that we have
\begin{equation*}
\int_{\Omega}\tn{S}(\nabla \vu):\nabla \vu_t\dx \geq C(\mu,\eta) \frac{d}{dt} \|\nabla \vu\|_{L^2}^2.
\end{equation*}
Moreover, the elliptic theory yields $\|\nabla^2 \vu\|_{L^2} \leq C \|A\vu\|_{L^2}$.
Now we apply the H\"older and
Young inequalities to most of the terms on the right hand side of \eqref{est_lin3} and
use the fact that $\vr$ is bounded from below by a positive constant. The first term
on the right hand side can be absorbed by the left hand side for sufficiently small $\ep$ and so we get
\begin{equation*}
\|\vu_t\|_{L^2}^2 + \frac{d}{dt}\|\nabla \vu\|_{L^2}^2 + \ep \|\nabla^2 \vu\|_{L^2}^2 \leq
C (\|\vc{V}_t\|_{H^1}^2 + \|\vc{\bar F}\|_{L^2}^2 + \|\vr \vu^b_{t}\|_{L^2}^2)
+ C\|\nabla \vu\|_{L^2}^2,
\end{equation*}
The energy estimate \eqref{ene} gives a bound on the last term of the right hand side. Therefore applying the Gronwall inequality we obtain
\begin{equation} \label{est_lin_2}
\|\vu_t\|_{L^2(L^2)}+\|\nabla \vu\|_{L^\infty(L^2)}+\|\nabla^2 \vu\|_{L^2(L^2)}\leq
\phi(\|\vr\|_{L^\infty(\Omega \times (0,T))})[\|\vc{V}_t\|_{L^2(H^1)} + \|\vc{F}\|_{L^2(L^2)} + \|\vu^b\|_{{\cal Y}(T)}].
\end{equation}
%
Next we take the time derivative of $\eqref{lin_ME1}$, multiply by $(\vu-\vc{V})_t$ and integrate. We get
\begin{align} \label{est_lin_4}
\frac{1}{2}\frac{d}{dt}\int_{\Omega}\vr|\vu_t|^2\dx + \int_{\Omega}A\vu_t\cdot (\vu-\vV)_t\dx =
&-\frac{1}{2}\int_{\Omega}\vr_t|\vu_t|^2\dx
+\int_{\Omega}\vr \vu_{tt}\cdot \vV_t\dx+\int_{\Omega}\vr_t \vu_t\cdot \vV_t\dx \nonumber\\
&+\int_{\Omega}(\vc{F}_t-\vr_t \vu^b_{t})\cdot (\vu-\vc{V})_t\dx
-\int_{\Omega}\vr \vu^b_{tt}\cdot (\vu-\vV)_t\dx.
\end{align}
As before we integrate by parts the second term on the left hand side.
Using the fact that $(\vu-\vc{V})_t\cdot \vn=0$ and the identity
$[\tn{S}(\nabla \vu_t)\vn]_{\rm tan}=\partial_t[\tn{S}(\nabla \vu)\vn]_{\rm tan} = 0$
we obtain
$$
\int_{\Omega}A\vu_t\cdot (\vu-\vV)_t\dx = \int_{\Omega}\tn{S}(\nabla \vu_t):\nabla \vu_t\dx
-\int_{\Omega}\tn{S}(\nabla \vu_t):\nabla \vV_t\dx .
$$
The condition $(\vu-\vc{V})_t\cdot \vn=0$ implies the Poincar\'e inequality for $(\vu-\vc{V})_t$
which yields
$$
\|\vu_t\|_{L^2} \leq C (\|\nabla \vu_t\|_{L^2}+\|\vc{V}_t\|_{H^1}).
$$
Now we examine the right hand side of \eqref{est_lin_4}.
For the first term we have by Poincar\'e inequality
$$
\left|\int_{\Omega}\vr_t|\vu_t|^2\dx\right| \leq \|\vr_t\|_{L^3}\|\vu_t\|_{L^6}(\|\nabla \vu_t\|_{L^2}+\|\vc{V}_t\|_{H^1}) \leq
\delta (\|\nabla \vu_t\|_{L^2}^2+\|\vc{V}_t\|_{H^1}^2) + C(\delta)\|\vr_t\|_{L^3}^2\int_{\Omega}\vr|\vu_t|^2\dx,
$$
where $\delta$ is a sufficiently small number coming from application of Young inequality
(we keep this notation in what follows).
The remaining terms are estimated directly and we get from \eqref{est_lin_4}
\begin{multline} \label{est_lin_5}
\frac{d}{dt}\int_{\Omega}\vr|\vu_t|^2\dx + \|\nabla \vu_t\|_{L^2}^2 \leq
C \Big\{ \delta_1 \|\nabla \vu_t\|_{L^2}^2 + C(\delta_1) \int_{\Omega}\vr|\vu_t|^2\dx
+\|\vV_t\|_{L^6}(\|\vr_t\|_{L^3}^2+\|\vu_t\|_{L^2}^2)  \\
+C(\delta_2)\phi(\|\vr\|_{L^\infty}) (\|\vc{\bar F}_t\|_{L^2}^2 + \|\vr_t \vu^b_{t}\|_{L^2}^2 +\|\vc{V}_t\|_{H^1}^2)
+\delta_2\|\vu_{tt}\|_{L^2}^2 \Big\}
-\int_{\Omega}\vr \vu^b_{tt}\cdot (\vu-\vV)_t\dx.
\end{multline}
The first term on the right hand side can be absorbed by the left hand side and the second term will be treated with the Gronwall
inequality. For the $L^2$ norm of $\vu_t$ we use \eqref{est_lin_2}.
Therefore integrating \eqref{est_lin_5} in time we obtain
\begin{align} \label{est_lin_6}
\|\vu_t\|_{L^{\infty}(0,T;L^2)}+\|\nabla \vu_t\|_{L^2(0,T;L^2)} \leq
&\phi(\|\vr\|_{L^\infty(H^2)},\|\vr_t\|_{L^2(H^1)})[\|\vV_t\|_{L^2(H^1)}+\phi_1(\|\vc{V}_t\|)_{L^4(L^2)}\|\vF_t\|_{L^2(L^2)}] \nonumber \\
&+\delta(\|\vu_{tt}\|_{L^2(L^2)}+\|\vV_{tt}\|_{L^2(L^2)})
+C(\delta)\|\vu^b_{tt} \vr\|_{L^2(L^2)} =: \Phi.
\end{align}
%
Next, computing $A\vu$ from \eqref{lin_ME1} we get a bound on $\|A\vu\|_{L^\infty(0,T;L^2)}$
which due to \eqref{est_lin_6} and ellipticity of $A$ yields
\begin{equation} \label{est_lin_9}
\|\nabla^2\vu\|_{L^\infty(L^2)} \leq
\|\vc{\bar F}\|_{L^\infty(L^2)} + \|\vr \vu^b_{t}\|_{L^\infty(L^2)}
+ \|\vr\|_{L^\infty((0,T)\times\Omega)}\Phi.
\end{equation}
%
%
In order to show the bound on $\|\vu\|_{L^2(H^3)}$
we take the spatial gradient of \eqref{lin_ME1} which yields
\begin{equation}
\nabla A\vu=\nabla \vc{\bar F} - \vr \nabla (\vu_t+\vu^b_{t}) - (\vu_t+\vu^b_{t}) \nabla \vr.
\end{equation}
We compute the $L^2(0,T;L^2)$ norm of the right hand side. From the H\"older inequality we get
$$
\|(\vu_t+\vu^b_{t}) \nabla \vr\|_{L^2(L^2)} + \|\vr\nabla(\vu_t+\vu^b_{t})\|_{L^2(L^2)}
\leq \phi(\|\vr\|_{L^\infty(H^2)}) \|\vu_t+\vu^b_{t}\|_{L^2(H^1)}.
$$
Therefore we obtain a bound on $\|\nabla^3 \vu\|_{L^2(L^2)}$ which,
together with previously obtained estimates on lower derivatives gives
a bound on $\|\vu\|_{L^2(H^3)}$ assuming we can estimate $\|\vu_{tt}\|_{L^2(L^2)}$.

Hence, in order to close the estimates we need to show a bound for $\|\vu_{tt}\|_{L^2(L^2)}$
which appears in the term $\Phi$ in \eqref{est_lin_6}.
Taking the time derivative of \eqref{lin_ME1} we obtain the elliptic problem for $\vu_t$:
\begin{align} \label{sys_w}
&-\div \tn{S}(\nabla \vu_t) =\vc{\bar F}_t-\vr_t\vu_t-\vr\vu_{tt}-\vr_t \vu^b_{t} - \vr\vu^b_{tt}, \nonumber\\[3pt]
&\left[ \tn{S}(\nabla \vu_t)\vc{n} \right]_{\rm tan}|_{\de \Omega} =0,\nonumber\\[3pt]
&(\vu-\vV)_t\cdot\vn|_{\de \Omega}=0,
\end{align}
recall that we have set $\kappa = 0$. For this problem the classical elliptic theory yields
\begin{equation} \label{ellip1}
\|\vu_t\|_{H^2} \leq C \left( \|\vc{\bar F}_t\|_{L^2}+\|\vr_t(\vu_t+\vu^b_{t})\|_{L^2}+\|\vr\vu_{tt}\|_{L^2}+\|\vr\vu^b_{tt}\|_{L^2} \right).
\end{equation}
Integrating \eqref{ellip1} with respect to time we get
\begin{equation} \label{e10}
\|\vu_t\|_{L^2(H^2)} \leq C \left( \|\vc{\bar F}_t\|_{L^2(L^2)}+\|\vu^b_{tt}\|_{L^2(L^2)}+\|\vu_{tt}\|_{L^2(L^2)}+\|\vr_t\|_{L^2(H^1)}(\|\vu_t\|_{L^\infty(H^1)}+\|\vu^b_{t}\|_{L^\infty(H^1)}) \right).
\end{equation}
Now we multiply \eqref{sys_w} by $(\vu-\vV)_{tt}$ and integrate over $\Omega$:
\begin{equation} \label{int_sys_w}
\int_{\Omega}\vr |\vu_{tt}|^2\dx-\int_{\Omega} \div \tn{S}(\nabla \vu_t)\cdot (\vu-\vV)_{tt}\dx =
\int_{\Omega} (\vc{\bar F}_t-\vr\vu^b_{tt})\cdot(\vu-\vV)_{tt}\dx
+ \int_{\Omega}\vr_t (\vu_t+\vu^b_{t}) \cdot(\vV-\vu)_{tt}
+ \vr \vu_{tt}\cdot \vV_{tt}\dx
\end{equation}
Again, from the boundary condition \eqref{sys_w}$_2$ we have
\begin{equation} \label{e100}
-\int_{\Omega} \div \tn{S}(\nabla \vu_t)\cdot (\vu-\vV)_{tt}\dx=
\int_{\Omega}\tn{S}(\nabla \vu_t):\nabla(\vu-\vV)_{tt}\dx
\end{equation}
and for the second term on the right hand side of \eqref{int_sys_w} we have
$$
\int_{\Omega}\vr_t (\vu+\vu^b)_t \cdot(\vV-\vu)_{tt}\dx
\leq \delta (\|\vu_{tt}\|_{L^2}^2+\|\vV_{tt}\|_{L^2}^2)
+ C(\delta)(\|\vr_t\|_{L^4}^2 \|\vu_t\|_{L^4}^2).
$$
Treating similarly the other terms on the right hand side of \eqref{int_sys_w} we obtain
\begin{equation}
\vr\|\vu_{tt}\|_{L^2}^2 + \frac{d}{dt}\|\nabla \vu_t\|_{L^2}^2 \leq \delta \|\vu_{tt}\|_{L^2}^2 +
C(\delta) \left(\|\vc{F}_t\|_{L^2}^2 + \|\vu^b_{tt}\|_{L^2}^2 + \|\vV_{tt}\|_{L^2}^2
+ \|\vr_t\|_{H^1}^2(\|\vu_t\|_{H^1}^2+\|\vu^b_{t}\|_{H^1}^2)\right).
\end{equation}
Integrating this inequality in time we get
\begin{equation} \label{e11}
\|\vu_{tt}\|_{L^2(L^2)} + \|\vu_t\|_{L^\infty(H^1)} \leq \phi(\|\vr_t\|_{L^2(H^1)}^2)
\left( \|\vF_t\|_{L^2(L^2)}+\|\vu^b_{tt}\|_{L^2(L^2)}+\|\vV_{tt}\|_{L^2(L^2)}+\|\vu_t(0)\|_{H^1}\right).
\end{equation}
\smallskip
\begin{Remark}
For $\kappa>0$ we obtain in \eqref{e100} additional boundary term
$$
\int_{\de \Omega}(\vu-\vV)_{tt}\cdot\kappa[(\vu-\vV)_t]_{\rm tan}\dS
$$
which contains $\vu_{tt}$.
However we can integrate this term in time
$$
\int_0^T\int_{\de \Omega} \frac{d}{dt}\kappa |\vu-\vV|_t^2\dS\dt = \int_{\de \Omega}\kappa |(\vu-\vV)_t(T)|^2\dS
-\int_{\de \Omega}\kappa |(\vu-\vV)_t(0)|^2\dS
$$
and the first term has good sign and second is given.
\end{Remark}

Finally we comment the term $\vu_t(0)$ appearing on the right hand side of \eqref{e11}.
To this end we differentiate \eqref{lin_ME1} with respect to space and multiply by
$\de_{x_i} \vu_t$ and formally take the resulting equation in $t=0$
(rigorously we show it for a smooth approximation and pass to the limit with the estimate)
obtaining
\begin{equation} \label{e12}
\|\vu_t(0)\|_{H^1} \leq C \left(\|\bar \vF(0)\|_{H^1}+\|\vu(0)\|_{H^3}\right).
\end{equation}
Combining \eqref{e10},\eqref{e11} and \eqref{e12} we get
\begin{equation*}
\|\vu_{tt}\|_{L^2(L^2)} + \|\vu_t\|_{L^\infty(H^1)} + \|\vu_t\|_{L^2(H^2)} \leq \phi(\|\vr_t\|_{L^2(H^1)}^2)
\left( \|\bar \vF_t\|_{L^2(L^2)}+\|\bar \vF(0)\|_{H^1}+\|\vu^b_{tt}\|_{L^2(L^2)}+\|\vV_{tt}\|_{L^2(L^2)}+\|\vu_t(0)\|_{H^1}\right),
\end{equation*}
which together with previous estimates and definition of $\vc{\bar F}$ gives \eqref{est_lin_ME_fix}.
Now, since \eqref{lin_ME1} is a linear parabolic problem,
the existence of solution follow from the estimates we have shown and classical theory of parabolic equations,
see for example \cite{V1}.

\qed

\subsection{Proof of Proposition \ref{Psol_mom}}
Let us denote
\begin{equation} \label{def_A}
A(t)=\|\vu\|_{{\cal{Y}}(t)}.
\end{equation}
We start with the estimate for the extension of the boundary data.
\bLemma{ext}
Let $\vV$ satisfy the assumptions of Theorem \ref{Tmain}.
Then there exists an extension $\vu^b(\vu,\vV)$ defined by \eqref{prob_ext} of the boundary data $\vc{d}(\vu,\vV)$, $\vc{B}(\vu,\vV)$ given by \eqref{bc1} and \eqref{bc2}
satisfying the estimate
\begin{equation} \label{est_ext}
\|\vu^b\|_{{\cal{Y}}(T)}\leq E(T)[1+\|(\vu-\vV)\|_{{\cal{Y}}(T)}],
\end{equation}
where $E(t)$ is continuous and $E(0)=0$.
\eL
The proof consist in defining the extension in a special way to ensure that, roughly speaking, the regularity of $\vu^b$ is the same as the regularity of $\vu$. We derive an explicit formula for $\vu^b$ and using the assumed regularity of $\vV$ and $\vu$ together with smallness of time we show the estimate \eqref{est_ext}. As the proof it is quite technical, we show it in the Appendix.

Next, we need the following estimate for the right hand side of the momentum equation in Lagrangian coordinates:
\bLemma{est_F}
For ${\bf R}(\vr,\vu)$ defined by \eqref{def_R} we have for any $\ep>0$
\begin{equation} \label{est_F}
\|\vr \vV \cdot \nabla_x \vu + {\bf R}(\vr,\vu)\|_{L^2(H^1)}
+\|\de_t[\vr \vV \cdot \nabla_x \vu + {\bf R}(\vr,\vu)]\|_{L^2(L^2)}
\leq \\
C[(\ep+\sqrt{t}C(\ep)+E(t))A(t)],
\end{equation}
where $E(t)$ is small for small times.
\eL
\bProof
The LHS of \eqref{est_F} contain derivatives w.r.t. $x$ and now we
need estimates in Lagrangian coordinates $y$. However we can easily observe that
$$
\vr \vV \cdot \nabla_x \vu =
\vr \vV \cdot \nabla_y \vu + R_1(\vr,\vu)
$$
where
$$
\|R_1(\vr,\vu)\|_{L^2(H^1)} + \|\de_t R_1(\vr,\vu)\|_{L^2(L^2)} \leq E(t)A(t)
$$
with $E(t)$ small for small times, therefore we can work with these terms directly.
Let us treat the term $\vr \vV \cdot \nabla \vu$.
From the interpolation inequality we have for any $\ep>0$
$$
\|\nabla^2 \vu\|_{L^2(L^4)} \leq \ep \|\vu\|_{L^2(H^3)} +  C(\ep)t\|\vu\|_{L^\infty(H^2)}.
$$
Therefore using again \eqref{est_lin_cont1} we get
$$
\int_0^t\int_{\Omega}|\vr|^2 |\vV|^2 |\nabla^2\vu|^2 \leq
\int_0^t \|\vr\|_{L^8}^2 \|\vV\|_{L^8}^2 \|\nabla^2 \vu\|_{L^4}^2 \leq
$$$$
C[\ep\|\vu\|_{L^2(H^3)}^2+tC(\ep)\|\vu\|_{L^\infty(H^2)}^2].
$$
Next, recalling the definition \eqref{def_R} of ${\bf R}(\cdot,\cdot)$ we see that it contains
terms with derivatives of $\vu$ of order up to $2$ multiplied by quantities which are small for small times, therefore
$$
\|{\bf R}(\vr,\vu)\|_{L^2(H^1)} \leq E(t)\|\vu\|_{L^2(H^3)},
$$
where $E(t)$ is small for small times.
Let us estimate the time derivative.
We have
$$
\int_0^t \int_{\Omega} \vr^2 |\vV|^2 |\nabla \vu_t|^2 \leq \|\vr\|_{L^\infty(\Omega \times (0,T))}\|\vV\|_{L^\infty(\Omega \times (0,T))}^2 t\|\vu_t\|_{L^\infty(H^1)}^2 \leq
Ct[A(t)]^2
$$
and treating similarly the other terms coming from the chain rule we obtain
\begin{equation*}
\|\de_t[\vr \vV \cdot \nabla \vu]\|_{L^2(L^2)}\leq C[\ep A(t)+C(\ep)tA(t)].
\end{equation*}
Finally, the structure of ${\bf R}$ implies
\begin{equation*}
\|\de_t{\bf R}(\vr,\vu)\|_{L^2(L^2)} \leq E(t)[\|\vu\|_{L^2(H^3)}+\|\vu_t\|_{L^2(H^2)}]\leq E(t)A(t)
\end{equation*}
with $E(t)$ as above. This estimate completes the proof of \eqref{est_F}.

\qed

Now for clarity let us denote again functions defined in Lagrangian coordinates by $\vut,\vrt$.
Combining \eqref{est_ext} and \eqref{est_F} we see that
the right hand side of the system \eqref{eq:mer2},\eqref{bc1},\eqref{bc2} satisfies
\begin{align}
&\|\tilde\vF(\vrt,\vut)\|_{L^2(H^1)}+\|\tilde\vF_t(\vrt,\vut)\|_{L^2(L^2)}
+\|\vu^b(\vut,\vV)\|_{{\cal{Y}}(T)} \leq
\|\vF\|_{L^2(H^1)}+\|\vF_t\|_{L^2(L^2)}\nonumber\\
&\quad+E(T)\left( A(t)+\|\vV\|_{L^2(H^3)}+\|\vV_t\|_{L^\infty(H^1)}+\|\vV_{tt}\|_{L^2(L^2)} \right)
\end{align}
where $\vu^b(\vut,\vV)$ is the extension of the boundary data $\vc{d}(\vut,\vV),\vc{B}(\vut,\vV)$ given by Lemma \ref{Lext}
and $E(t)$ is small for small times. Therefore from \eqref{est_lin_ME_fix} we obtain the estimate
\eqref{est_lin_ME}. Moreover, the right hand side of \eqref{eq:mer2}, \eqref{bc1}, \eqref{bc2} is linear w.r.t. $\vut$,
therefore the existence of a unique solution follows from the Banach fixed point theorem. Therefore, as Lagrangian transformation is a diffeomorphism for small times, $\vu(t,x)=\vut(t,Y(y,x))$ is a solution of \eqref{lin_ME}.

\section{Proof of Theorem \ref{Tmain}}\label{s:proof}

\subsection{Boundedness of the sequence of approximations}
In this section we use the estimates for the linear problems to show that the sequence
$(\vr_n,\vu_n)$ defined by the iterative scheme described in Section \ref{ss:scheme} is bounded in the space where we are looking for the solution.
Let us denote
\begin{equation} \label{def_An}
A_n(t)=\|\vu_n\|_{{\cal{X}}(T)}.
\end{equation}
Since our estimate for the linear momentum equation holds on fixed domain, we rewrite
\eqref{eq:MEn} in Lagrangian coordinates
\begin{align}\label{eq:MEnf}
\vr_{n+1} \de_t \vu_{n+1}
- \mu\Delta_y \vu_{n+1} - (\frac{\mu}{3} + \eta)\nabla_y{\rm div}_y\vu_{n+1}&={\bf R}_n,\nonumber\\[3pt]
(\vu_{n+1}-\vc{V})\cdot\vc{n}|_{\Gamma}&=\vc{d}(\vu_n,\vV),\nonumber\\[3pt]
\left([\tn{S}(\nabla_x \vu_{n+1})\vn]\cdot \tau^k + \kappa [\vu_{n+1}-\vV]\cdot \tau^k\right)|_{\Gamma}&=\vc{B}(\vu_n,\vV),
\end{align}
where
\begin{equation} \label{def_Rn}
{\bf R}_n=\vr_{n+1} (\vV \cdot \nabla_y \vu_{n+1} - \vu_n \cdot \nabla_x \vu_{n}) - \nabla_x p(\vr_{n+1}) + {\bf R}(\vr_{n+1},\vu_{n+1})
\end{equation}
and ${\bf R}(\cdot,\cdot)$ is defined in \eqref{def_R}.
The following lemma gives the estimate $L^2(H^1)$ norm of the right hand side of \eqref{eq:MEnf}$_1$ and $L^2(L^2)$ norm of its time derivative.
\bLemma{est_F1}
For ${\bf R}(\vr_{n+1},\vu_{n+1})$ and ${\bf R}_n$ defined by \eqref{def_R} and \eqref{def_Rn} respectively we have for any $\ep>0$
\begin{equation} \label{est_Fn}
\|{\bf R}_n\|_{L^2(H^1)}
+\|\de_t{\bf R}_n\|_{L^2(L^2)}
\leq
\phi(\sqrt{t}A_n(t))[(\ep+\sqrt{t}C(\ep)+E(t))A_n(t)+E(t)A_{n+1}(t)],
\end{equation}
where $E(t)$ is small for small times.
\eL
\bProof
Arguing as in the proof of \eqref{est_F} we can replace the derivatives w.r.t. $x$
in the definition of ${\bf R}_n$ with derivatives w.r.t. $y$.
Since the pressure $p$ is a $C^2$ function of the density we have
$$
\nabla^2 p(\vr)=p''(\vr)|\nabla \vr|^2+p'(\vr)\nabla^2\vr \sim \vr (|\nabla \vr|^2+\nabla^2 \vr),
$$
and therefore
\begin{equation}
\|\nabla p(\vr)\|_{L^2(H^1)}^2 \leq \|\vr\|_{L^\infty(L^2)}^2 t \|\vr\|_{L^\infty(H^2)}^2(1+\|\vr\|_{L^\infty(H^2)}^2) \leq
\phi(\sqrt{t}\|\vu\|_{L^2(H^3)})t \|\vr\|_{L^\infty(L^2)}^2,
\end{equation}
where in the last passage we applied \eqref{est_lin_cont1}.
Next, similarly to the proof of \eqref{est_F} we get
\begin{equation}
\|\vr_{n+1} (\vu_n \cdot \nabla \vu_n - \vV \cdot \nabla \vu_{n+1}) \|_{L^2(H^1)} \leq
\phi(\sqrt{t}\|\vu\|_{L^2(H^3)})
[\ep\|\vu\|_{L^2(H^3)}+tC(\ep)\|\vu\|_{L^\infty(L^2)}].
\end{equation}
and
$$
\|{\bf R}(\vr_{n+1},\vu_{n+1})\|_{L^2(H^1)} \leq E(t)\|\vu_{n+1}\|_{L^2(H^3)},
$$
where $E(t)$ is small for small times.
Let us estimate the time derivative. For the pressure we have
$$
\de_t \nabla p(\vr)\sim \vr(\nabla \vr_t + \vr_t \nabla \vr),
$$
therefore from \eqref{est_lin_cont1} and \eqref{est_lin_cont2} we obtain
\begin{equation}
\|\de_t \nabla p(\vr)\|_{L^2(L^2)}\leq \phi(\sqrt{t}A_n(t))[\ep A_n(t)+C(\ep)tA_n(t)].
\end{equation}
The remaining terms are again estimated like in the proof of \eqref{est_F}
which leads to
\begin{equation}
\|\de_t[\vr_n (\vu_n \cdot \nabla \vu_n - \vV \cdot \nabla \vu_{n+1})]\|_{L^2(L^2)}\leq \phi(\sqrt{t}A_n(t))[\ep A_n(t)+C(\ep)tA_n(t)].
\end{equation}
and
\begin{equation}
\|\de_t{\bf R}(\vr_{n+1},\vu_{n+1})\|_{L^2(H^1)} \leq E(t)[\|\vu_{n+1}\|_{L^2(H^3)}+\|\vu_{n+1,t}\|_{L^2(H^2)}]\leq E(t)A_{n+1}(t)
\end{equation}
with $E(t)$ as above. This estimate completes the proof of \eqref{est_Fn}.

\qed

With above lemmas we are ready to show the key estimate for the sequence of approximations
\bProposition{est_An}
Let $A_n(t)$ be defined in \eqref{def_An}. Then there exists $M>0$ sufficiently large
and $T^*>0$ such that
\begin{equation} \label{est_an}
A_n(t) \leq M \quad \textrm{for} \quad t \leq T^* .
\end{equation}
\eP
\bProof
The appropriate choice of the extension of the boundary data given by Lemma \ref{Lext} ensures
$$
\|\vu^b_{tt}(\vu_n,\vV)\|_{{\cal{Y}}(T)} \leq E(T) A_n(t).
$$
Therefore the estimate \eqref{est_lin_ME} applied to \eqref{eq:MEnf} yields
$$
A_{n+1}(t) \leq \phi(\sqrt{t}A_n(t))\big[(\ep+tC(\ep)+E(t))A_n(t)
+\|\vu_0\|_{H^3}+\|\vr_0\|_{H^2}+\|\vV\|_{{\cal Y}(T)}\big]
$$
for some increasing positive function $\phi(\cdot)$. We conclude that there exists $M=M(\vu_0,\vV)$ and
$T^*>0$ such that \eqref{est_an} holds. \qed

\subsection{Convergence of the sequence of approximations.}
Let us denote
$$
\vw_{n+1} = \vu_{n+1} - \vu_n, \quad \sigma_{n+1} = \vr_{n+1} - \vr_n.
$$
Subtracting \eqref{eq:MEn} for $n+1$ and $n$ we get
\begin{align} \label{eq:MEdif}
&\vr_{n}\de_t \vw_{n+1} - \mu \Delta \vw_{n+1}-(\frac{\mu}{3}+\eta)\nabla\ddiv \vw_{n+1}= \nonumber\\
&=\sigma_{n+1}\de_t \vu_{n+1}+\sigma_{n+1}(\vu_{n-1}\cdot\nabla\vu_{n-1})
+\vr_{n+1}(\vu_n \cdot \nabla \vw_n + \vw_n\cdot\nabla\vu_{n-1} )\\ \nonumber
&+p'(\vr_n)\nabla \sigma_n+\nabla \vr_{n-1}\sigma_n\int_0^1 p''(\vr_n)(s\vr_n+(1-s)\vr_{n-1})ds
\end{align}
supplied with boundary conditions
\begin{equation}
\vw_{n+1}\cdot\vc{n}|_{\Gamma_t}=0, \qquad
[\tn{S}(\nabla \vw_{n+1})\vn]_{\rm tan} + \kappa [\vw_{n+1}]_{\rm tan}|_{\Gamma_t}=0.
\end{equation}
Subtracting \eqref{eq:CEn} we obtain
\begin{equation} \label{eq:CEdif}
\de_t \sigma_{n+1}+\vu_n\cdot\nabla \sigma_{n+1}+\sigma_{n+1}\ddiv \vu_n=
-\vr_n \ddiv \sigma_n-\sigma_n\cdot \nabla\vr_n, \quad \sigma_{n+1}(0,\cdot)=0.
\end{equation}
In order to apply our estimates for the linear momentum equation
we rewrite \eqref{eq:MEdif} on a fixed domain using the transformation \eqref{eq:coc}.
We obtain
\begin{align} \label{eq:MEdif_lag}
&\vr_{n}\de_t \vw_{n+1} - \mu \Delta_y \vw_{n+1}-(\frac{\mu}{3}+\eta)\nabla\ddiv_y \vw_{n+1}=
\vr_n \vV \cdot \nabla_y \vw_{n+1} + {\bf R}(\vr_n,\vw_{n+1})   \nonumber\\
&+\sigma_{n+1}\de_t \vu_{n+1}+\sigma_{n+1}(\vu_{n-1}\cdot\nabla_x\vu_{n-1})
+\vr_{n+1}(\vu_n \cdot \nabla_x \vw_n + \vw_n\cdot\nabla_x\vu_{n-1} )\\ \nonumber
&+p'(\vr_n)\nabla_x \sigma_n+\nabla \vr_{n-1}\sigma_n\int_0^1 p''(\vr_n)(s\vr_n+(1-s)\vr_{n-1})ds
=: {\bf{\tilde R}}_n
\end{align}
with boundary conditions
\begin{align} \label{bc_dif}
&(\vw_{n+1} \cdot \vc{n})(y)|_{\Gamma} = \vw_{n+1} \cdot\left(\vc{n}(y)-\vc{n}({X}(t,y))\right)=:\vc{d_0}(\vw_{n+1}), \\ \nonumber
&\mu \vn(y)(\nabla_y \vu+\nabla^t_t \vu) \cdot \tau^k(y) = \mu\left(\nabla_y\vw_{n+1}(t,y)(\tn{I}-\nabla_x\vc{Y}) + ((\tn{I}-\nabla_x^T\vc{Y})\nabla_y^T\vw_{n+1}(t,y))^T\right)\vc{n}(\vc{X}(t,y))\cdot\tau^k(\vc{X}(t,y)) \\ \nonumber
&\qquad + \mu(\nabla_y\vw_{n+1}+\nabla^T_y\vw_{n+1})(t,y)[(\vc{n}(y)-\vc{n}(\vc{X}(t,y)))\cdot\tau^k(\vc{X}(t,y)) + \vc{n}(t,y)\cdot(\tau^k(y)-\tau^k(\vc{X}(t,y)))] \\ \nonumber
&\qquad + \kappa \vw_{n+1}(t,y)\cdot(\tau^k(y)-\tau^k(\vc{X}(t,y)))  =: \vc{B_0}(\vw_{n+1})(t,y).
\end{align}
The left hand side of the system \eqref{eq:MEdif_lag} has exactly the structure of \eqref{lin_ME_fix}.
As we will not be able to close the estimate for ${\bf w}_n$ in the regularity we have for the sequence ${\bf u}_n$ (see Remark \ref{r:conv} below) we show the convergence in a weaker space. Let us denote
\begin{equation}\label{norm:bn}
B_n(t)=\|\vw_n\|_{L^\infty(H^1)}+\|\vw_n\|_{L^2(H^2)}+\|\vw_{n,t}\|_{L^2(L^2)}.
\end{equation}
Repeating the proof of \eqref{est_lin_ME} we obtain in particular the following (in fact classical)
parabolic estimate
\begin{equation} \label{est_Bn}
B_n(t) \leq C [\|{\bf{\tilde R}}_n\|_{L^2(L^2)}+
\|\vc{B_0}(\vw_{n+1})\|_{L^2(H^{1/2}(\de \Omega))}+\|{\bf d_0}(\vw_{n+1})\|_{L^2(H^{3/2}(\de \Omega))}].
\end{equation}
The structure of the boundary data clearly implies
\begin{equation}
\|\vc{b}(\vw_{n+1})\|_{L^2(H^{1/2}(\de \Omega))}+\|{\bf d}(\vw_{n+1})\|_{L^2(H^{3/2}(\de \Omega))}\leq
E(t) \|\vw_{n+1}\|_{L^2(H^2)}.
\end{equation}
The following lemma gives the estimate of the right hand side of \eqref{eq:MEdif_lag}.
\bLemma{est_Rdif}
Let ${\bf{\tilde R}}_n$ be defined in \eqref{eq:MEdif_lag}. Then
\begin{equation} \label{est_Rdif}
\|{\bf{\tilde R}}_n\|_{L^2(L^2)}\leq C[\|\sigma_{n+1}\|_{L^\infty(L^2)}+\|\sigma_n\|_{L^\infty(L^2)}
+(\ep+tC(\ep)+E(t))(\|\nabla \sigma_n\|_{L^\infty(L^2)}+B_n(t))],
\end{equation}
where $E(t)$ is small for small times.
\eL
\bProof
Arguing as in the proofs of Lemma \ref{Lest_F1} and Proposition \ref{Pest_An} we get
\begin{equation}
\|{\bf R}(\vr_n,\vw_{n+1})\|_{L^2(H^1)} \leq E(t)\|\vw_{n+1}\|_{L^2(H^3)}.
\end{equation}
For the remaining terms, we can follow the proof of Lemma \ref{Lest_F1}.
As we have explained there, we can replace the derivatives w.r.t $x$ with derivatives w.r.t $y$.
We have
\begin{align}
&\int_0^t \int_{\Omega} |\sigma_{n+1}|^2 |\vu_{n-1}|^2 |\nabla \vu_{n-1}|^2 \leq
\|\sigma_{n+1}\|_{L^\infty(L^2)}^2 \|\vu_{n-1}\|_{L^\infty(H^2)}^2\int_0^t \|\nabla \vu_{n-1}\|_{L^\infty}^2 \leq \\ \nonumber
& 
\|\sigma_{n+1}\|_{L^\infty(L^2)}^2 \|\vu_{n-1}\|_{L^\infty(H^2)}^2 \left[ \ep\|\vu_{n-1}\|_{L^2(H^3)}^2+C(\ep)t\|\vu_{n-1}\|_{L_\infty(H^2)}^2 \right] \leq  \nonumber\\
&C \| \sigma_{n+1}\|^2_{L^\infty(L^2)}[(\ep+tC(\ep))A_{n-1}(t)]^2.
\end{align}
Next, consider the term
\begin{equation}
\|\vr_n \vV \cdot \nabla_y \vw_{n+1}\|_{L^2(L^2)}^2 \leq
\|\vr_n\|_{L^\infty(H^2)}^2 \|\vV\|_{L^\infty(H^2)}^2 \int_0^t \|\nabla \vw_{n+1}\|_{L^2}^2\leq
C[ \ep \|\vw_{n+1}\|_{L^2(H^2)} + tC(\ep)\|\vw_{n+1}\|_{L^\infty(L^2)} ].
\end{equation}
Similarly
\begin{equation}
\|\vr_{n+1}(\vu_n \cdot \nabla \vw_n+\vw_n\cdot \nabla \vu_{n-1}) \|_{L^2(L^2)}\leq C[ \ep \|\vw_{n+1}\|_{L^2(H^2)} + tC(\ep)\|\vw_{n+1}\|_{L^\infty(L^2)} ]
\end{equation}
and
\begin{equation}
\|\sigma_{n+1} \de_t \vu_{n+1}\|_{L^2(L^2)} \leq \|\sigma_{n+1}\|_{L^\infty(L^2)}\|\de_t \vu_{n+1}\|_{L^2(H^2)}.
\end{equation}
Finally, under our regularity assumptions on the pressure we have
\begin{equation}
\|p'(\vr_n)\nabla \sigma_n\|_{L^2(L^2)}^2 \leq C \|\nabla \sigma_n\|^2_{L^\infty(L^2)}\int_0^t\|p'(\vr_n)\|_{L^2}^2\leq
C \|\sigma_n\|^2_{L^\infty(H^1)}[\ep\|\vr_n\|^2_{L^2(H^2)}+tC(\ep)\|\vr\|^2_{L^\infty(L^2)}]
\end{equation}
and
\begin{equation}
\|\nabla \vr_{n-1}\sigma_n\int_0^1 p^{''}(\vr_n)(s\vr_n+(1-s)\vr_{n-1})ds\|_{L^2(L^2)}\leq
C\|\sigma_n\|_{L^\infty(L^2)}\|\nabla \vr_{n-1}\|_{L^2(L^\infty)}.
\end{equation}
Combining all above estimates and Proposition \ref{Pest_An} we get \eqref{est_Rdif}.

\qed

\noindent
In order to show the convergence we need to estimate the norms of $\sigma_n$ on the right hand side of \eqref{est_Rdif}.
For this purpose we investigate the equation \eqref{eq:CEdif}.
Again, as in the proof of Proposition \ref{Psol_cont}, we stay in a moving domain and
use the transformation
\begin{equation} \label{char_dif}
X(t,z)=z+\int_0^t \vu_n(X(s,z)){\rm d}s.
\end{equation}
Denoting
$$
\bar \sigma_{n+1}(t,z)=\sigma_{n+1}(t,X(t,z)),
$$
equation \eqref{eq:CEdif} rewrites as a nonhomogeneous ODE
\begin{equation} \label{CEdif_lag}
\de_t \bar \sigma_{n+1}+\bar \sigma_{n+1} \ddiv_x \vu_n = - (\vr_n \ddiv \vw_n+\vw_n\cdot \nabla\vr_n),
\quad \bar \sigma_{n+1}(0,\cdot)=0.
\end{equation}
The solution for this equation is given by the following explicit formula
\begin{align} \label{sol_form}
&\bar \sigma_{n+1}(t,z)={\rm exp}[-\int_0^t \ddiv_x \vu_n(s,X(x,z))ds]\cdot\\ \nonumber
&\int_0^t \left\{ - \Big(\vr_n \ddiv_x \vw_n+\vw_n\cdot \nabla_x\vr_n)(r,X(r,z)\Big)
{\rm exp}[-\int_0^r \ddiv_x \vu_n(\tau,X(\tau,z))d\tau] \right\}dr.
\end{align}
\begin{Remark}\label{r:conv}
Notice that the above formula contains $\nabla_x \vr_n$ and this is the reason why we have to work with lower regularity. Namely, if we would like to have the regularity ${\cal X}(T)$
then we would have to estimate $L^2(L^2)$ of the derivatives of the right hand side of \eqref{eq:MEdif_lag}, therefore $\nabla_x\sigma_n$ which would require information about $\nabla_x^3 \vr_n$. Such situation when lack of sufficient information on the density requires to work in lower regularity is typical for the compressible Navier-Stokes system, see for example \cite{Za},\cite{PP}.
\end{Remark}
The required estimate is provided in the following
\bLemma{est_CEdif}
Let $\sigma_{n+1}$ solve \eqref{eq:CEdif}. Then
\begin{equation} \label{est_CEdif}
\|\sigma_{n+1}\|_{L^\infty(L^2)} \leq (\ep + tC(\ep))B_n(t), \qquad
\|\nabla \sigma_{n+1}\|_{L^\infty(L^2)} \leq C B_n(t).
\end{equation}
\eL
\bProof The proof relies on the formula \eqref{sol_form} and follows the proof of Proposition \ref{Psol_cont}, therefore here we show only the main ideas.
The solution formula \eqref{sol_form} clearly implies
\begin{align}
&\|\bar \sigma_n(t,z)\|_{L^2} \leq C \int_0^t \|(\vr_n \ddiv \vw_n)(s,\cdot)\|_{L^2}+\|(\vw_n \cdot \nabla \vr_n)(s,\cdot)\|_{L^2} ds \leq \\ \nonumber
&C\|\vr_n\|_{L^\infty(H^1)}\int_0^t \|\vw_n\|_{H^1}ds \leq C \int_0^t \ep \|\nabla^2 \vw_n\|_{L^2}+C(\ep)\|\vw_n\|_{L^2}\leq
C [\ep B_n(t)+tC(\ep)B_n(t)],
\end{align}
and the regularity of the change of coordinates together with smallness of time yields the
same estimate for $\sigma_n(t,x)$ which gives the first statement of \eqref{est_CEdif}.
Let us denote
\begin{equation}
\ve(t,z)={\rm exp}[-\int_0^t \ddiv_x \vu_n(s,X(s,z))ds].
\end{equation}
Differentiating \eqref{sol_form} w.r.t $z$ we get
$$
\nabla_x \sigma_{n+1} = \nabla_x z \left\{
\underbrace{\nabla_z \left[ \ve(t,z) \right] \int_0^t \left\{ - \Big(\vr_n \ddiv_x \vw_n+\vw_n\cdot \nabla_x\vr_n)(r,X(r,z)\Big)
\ve(r,z) \right\}dr}_{I_1} \right.
$$$$
+ \left.\underbrace{\ve(t,z) \cdot \nabla \left[ \int_0^t \left\{ - \Big(\vr_n \ddiv_x \vw_n+\vw_n\cdot \nabla_x\vr_n)(r,X(r,z)\Big)
\ve(r,z) \right\}dr  \right]}_{I_2}\right\}.
$$
Applying the estimate \eqref{nabla_x_1} to the change of coordinates \eqref{char_dif}
we obtain
$$
\|\nabla_x z\,I_1\|_{L^2} \leq C \|[-\int_0^t \ddiv_x \vu_n(s,X(s,z))ds]\|_{L^\infty}^2 \|\int_0^t \nabla \ddiv \vu_n\|_{L^4} \|\int_0^t \{-\vr_n \ddiv \vw_n + \vw_n \cdot \nabla \vr_n\}ds\|_{L^4},
$$
where all the derivatives are already w.r.t. $x$. Therefore
\begin{align} \label{I1}
&\|\nabla_x z\,I_1\|_{L^2} \leq C \left(\int_0^t \|\nabla \ddiv \vu_n\|_{L^4}ds\right) \int_0^t \{\|\vr_n \ddiv \vw_n\|_{L^4}+\|\vw_n \cdot \nabla \vr_n\|_{L^4}\}ds \leq \nonumber \\
& C[ \ep \|\vu_n\|_{L^2(H^3)}+C(\ep)\|\vu_{n+1}\|_{L^\infty(L_2)} ]\|\vr_n\|_{L^\infty(\Omega \times (0,T))}\|\vw_n\|_{L^2(H^2)} \leq
C [\ep + t C(\ep)] B_n(t).
\end{align}
With $I_2$ we have
$$
I_2 = \underbrace{\ve(t,z)\int_0^t \left\{(-\int_0^r \nabla_z \ddiv_x \vu_n)\ve(r,z)(-\vr_n \ddiv \vw_n+\vw_n \cdot \nabla \vr_n)\right\}dr}_{I_{21}}
$$$$
+ \underbrace{\ve(t,z)\int_0^t \left\{ \ve(r,z)(-\nabla \vr_n \ddiv \vw_n-\vr_n\nabla\ddiv \vw_n + \nabla\vw_n\cdot\nabla\vr_n+\vw_n \cdot \nabla^2 \vr_n)\right\}dr}_{I_{22}}.
$$
The first part can be estimated exactly as $I_1$:
\begin{equation} \label{I21}
\|\nabla_x z\,I_{21}\|_{L^2} \leq C [\ep + t C(\ep)] B_n(t).
\end{equation}
However, in $I_{22}$ we cannot get smallness in time because of the presence of second derivatives of $\vw_n$,
namely
$$
\int_0^t \|\vr_n \nabla \ddiv \vw_n\|_{L^2} \leq \|\vr_n\|_{L^\infty(\Omega \times (0,T))} \|\nabla^2 \vw_n\|_{L^1(L^2)}\leq
C B_n(t),
$$
and therefore
\begin{equation} \label{I22}
\|\nabla_x z\,I_{22}\|_{L^2} \leq C B_n(t).
\end{equation}
Combining \eqref{I1}, \eqref{I21} and \eqref{I22} we obtain the second statement of \eqref{est_CEdif}.

\qed

We are now in a position to prove the following Lemma which gives the Cauchy condition
for the sequence $B_n(t)$, and therefore completes the proof of convergence of the sequence of approximations.
\bLemma{est_bn}
There exists $T>0$ and $0<K<1$ such that
\begin{equation} \label{est_Bn_final}
B_n(t) \leq KB_{n-1}(t) \quad {\rm for} \quad t\leq T^*.
\end{equation}
\eL
\bProof Combinining \eqref{est_Bn}, \eqref{est_Rdif} and \eqref{est_CEdif} we obtain
$$
B_n(t) \leq (\ep + t(\ep) + E(t))(B_n(t)+B_{n-1}(t))
$$
which implies \eqref{est_Bn_final}.

\qed

Now we conclude the proof of Theorem \ref{Tmain} in a standard way. The sequence $\vu_n$ converges strongly in the space with topology given by the norm \eqref{norm:bn}. On the
other hand, the bound \eqref{est_an} implies weak convergence up to a subsequence in ${\cal X}(T)$.
The estimates for the linear transport equations imply analogous convergence for the sequence of densities in appropriate spaces.
Therefore the limit $(\vu,\vr)$ satisfies the regularity given in Theorem \ref{Tmain}.

\section{Weak-strong uniqueness} \label{s:WSU}
\subsection{Energy inequality}\label{ei}
In the remainder of the paper we assume that $\kappa,\eta = 0$, however all the results hold (with appropriate modifications of formulae) also with $\kappa,\eta > 0$.
Our first result is a crucial observation allowing for proving the main theorem about weak-strong uniqueness.

\begin{Proposition}\label{p:ex2}
Let the assumptions of Theorem \ref{t:ex} be satisfied. Then the problem \eqref{i1a}-\eqref{i6c}, \eqref{ic:weak} admits a weak solution on any time interval $(0,T)$ in the sense specified through Definition \ref{d:ws}, satisfying moreover the energy inequality in the following form
\begin{align}\label{eq:EI}
&\int_{\Omega_\tau} \left( \frac{1}{2} \vr |\vu|^2 + H(\vr) \right)(\tau, \cdot) \ \dx +
\frac{1}{2} \int_0^\tau \int_{\Om_t} \mu \left| \Grad \vu + \Grad^t \vu - \frac{2}{3} \Div \vu \tn{I} \right|^2 \ \dxdt \\ \nonumber
\leq &\int_{\Omega_0} \left( \frac{1}{2 \vr_0} |(\vr \vu)_0 |^2 + H(\vr_0) \right) \ \dx + \int_{\Om_\tau} (\vr \vu \cdot \vc{V}) (\tau, \cdot) \ \dx - \int_{\Om_0} (\vr \vu)_0 \cdot \vc{V}(0, \cdot) \ \dx \\ \nonumber
+ &\int_0^\tau \int_{\Om_t} \left( \mu \left(\Grad \vu + \Grad^t \vu - \frac{2}{3} \Div \vu \tn{I} \right) : \Grad \vc{V} - \vr \vu \cdot \partial_t \vc{V} - \vr \vu \otimes \vu : \Grad \vc{V} - p(\vr) \Div \vc{V} \right)  \dxdt.
\end{align}
\end{Proposition}

\bProof
We follow the same series of approximations and penalizations as it is introduced in \cite[Section 3]{FKNNS} in the proof of Theorem \ref{t:ex}, see also  the Appendix for more details. The starting point is thus the modified energy inequality written on the fixed domain $B$ which is a ball large enough such that $\Omega_t \subset B$ for all $t \in [0,T]$ and $\vc{V} = 0$ on $\partial B$, see formula (3.10) in \cite{FKNNS}
\bFormula{p7}
\int_B \left( \frac{1}{2} \vr |\vu|^2 + H(\vr) + \frac{\delta}{\beta - 1} \vr^\beta \right)(\tau, \cdot) \ \dx +
\frac{1}{2} \int_0^\tau \int_B \mu_\omega \left| \Grad \vu + \Grad^t \vu - \frac{2}{3} \Div \vu \tn{I} \right|^2 \ \dxdt
\eF
\[
+
\frac{1}{\ep} \int_0^\tau \int_{\Gamma_t} \left| \left( \vu - \vc{V} \right) \cdot \vc{n} \right|^2 \ {\rm dS}_x \ \dt
\leq \int_B \left( \frac{1}{2 \vr_{0,\delta}} |(\vr \vu)_{0,\delta} |^2 + H(\vr_{0, \delta}) + \frac{\delta}{\beta - 1} \vr_{0, \delta}^\beta \right) \ \dx
\]
\[
+ \int_B \Big(  (\vr \vu \cdot \vc{V}) (\tau, \cdot) - (\vr \vu)_{0,\delta} \cdot \vc{V}(0, \cdot) \Big)
\ \dx
\]
\[
+ \int_0^\tau \int_B \left( \mu_\omega \left(\Grad \vu + \Grad^t \vu - \frac{2}{3} \Div \vu \tn{I} \right) : \Grad \vc{V} - \vr \vu \cdot \partial_t \vc{V} - \vr \vu \otimes \vu : \Grad \vc{V} -
p(\vr) \Div \vc{V} - \frac{ \delta }{\beta - 1} \vr^\beta \Div \vc{V} \right)  \dxdt.
\]
Let us remark that in \eqref{p7} we use the identity
$$ \tn S(\nabla_x \vu) : \nabla_x \vu = \left[\mu_\omega \left( \Grad \vu + \Grad^t \vu\right) - \frac{2}{3}\mu_{\omega} \Div \vu \tn{I} \right] : \nabla_x \vu = \frac 12 \mu_\omega \left| \Grad \vu + \Grad^t \vu - \frac{2}{3} \Div \vu \tn{I} \right|^2.  $$
Passing first with $\ep$ to zero, it is not difficult to observe that using the a priori estimates available, all the terms on the right hand side of \eqref{p7} converge to their counterparts. On the left hand side the last term is positive and thus can be omitted. Finally, using the convexity of $\tn{S}_\omega(\Grad\vu):\Grad\vu$ we have
\begin{equation}\label{eq:convex}
	\int_0^T \int_{B} \tn{S}_\omega( \Grad \vu) :\Grad \vu \, \dxdt \leq \liminf\limits_{\ep \to 0}  \int_0^T \int_{B} \tn{S}_\omega(\Grad \vu_\ep) :\Grad \vu_\ep \, \dxdt.
	\end{equation}

Next, passing with $\omega$ to zero, we first observe that all the terms which include the density can be rewritten as the integrals over $\Omega_t$ instead of integrals over $B$ using the fundamental Lemma 4.1 in \cite{FKNNS} due to artificial pressure term $\delta \rho ^{\beta}$, $\beta >4$, which gives us better integrability of the density. The viscosity term on the right hand side can be treated easily, in particular the integral over $B \setminus \Omega_t$ vanish due to the fact that $\mu_\omega \sil 0$ on this set. On the left hand side we split the integral of the term with stress tensor into two parts, the integral over $B \setminus \Omega_t$ can be omitted since it is positive and on $\Omega_t$ we use the fact that $\mu_\omega = \mu$ is constant and thus we can use again the convexity of $\tn{S}(\Grad\vu):\Grad\vu$ to obtain similar inequality as \eqref{eq:convex}.

Finally, we pass with $\delta$ to zero. Here we use the - nowadays already standard - results of \cite{EF70} to pass to the limit on the right hand side and to adjust the initial conditions, while on the left hand side we use the weak lower semicontinuity of the energy at time $\tau$.

\qed

\subsection{Relative energy inequality}\label{rei}

Having now already the energy inequality, we can deduce the relative energy inequality in the spirit of \cite{FeJiNo}. Before stating the theorem, we introduce some notation. For a weak solution $(\vr,\vu)$ and a pair of test functions $(r,\vU)$ defined on $Q_T$ we define the relative energy $\Ecal\Big([\vr,\vu]|[r,\vU]\Big)$  as
\begin{equation}\label{eq:RE}
\Ecal\Big([\vr,\vu]|[r,\vU]\Big)(\tau) = \int_{\Omega_{\tau}} \left( \frac 12 \vr\abs{\vu-\vU}^2 + H(\vr) - H'(r)(\vr-r)-H(r) \right)(\tau,\cdot)\, \dx.
\end{equation}

We prove the following

\begin{Proposition}\label{p:REI}

Let $(\vr,\vu)$ be a weak solution to the compressible Navier-Stokes system \eqref{i1a}-\eqref{i6c}, \eqref{ic:weak} constructed in Proposition \ref{p:ex2}. Then $(\vr,\vu)$ satisfies the following relative energy inequality
\begin{align}\label{eq:REI}
&\Ecal\Big([\vr,\vu]|[r,\vU]\Big)(\tau) + \int_0^\tau\int_{\Om_t} \left(\tn{S}(\Grad \vu) - \tn{S}(\Grad \vU)\right):\left(\Grad \vu - \Grad \vU\right)\,\dxdt \\ \nonumber \leq\, &\Ecal\Big([\vr_0,\vu_0]|[r(0,\cdot),\vU(0,\cdot)]\Big) + \int_0^\tau \Rcal(\vr,\vu,r,\vU)(t) \dt
\end{align}
for a.a. $\tau \in (0,T)$ and any pair of test functions $(r,\vU)$ such that $\vU \in C^{\infty}_c(\overline{Q_T})$, $\vU\cdot\n = \vV\cdot\n$ on $\Gamma_t$ for $t \in [0,T]$, $r \in C^{\infty}_c(\overline{Q_T})$, $r > 0$.
The remainder term $\Rcal$ is given by
\begin{align}\label{eq:R}
\Rcal(\vr,\vu,r,\vU)(t) &= \int_{\Om_t} \vr(\partial_t\vU + \vu\cdot\Grad\vU)\cdot(\vU-\vu) + \tn{S}(\Grad\vU):(\Grad\vU-\Grad\vu) \,\dx \\ \nonumber
&+ \int_{\Omega_t} \Div\vU(p(r)-p(\vr)) + (r-\vr)\partial_tH'(r) + (r\vU-\vr\vu)\cdot\Grad H'(r)\, \dx
\end{align}
\end{Proposition}
\bProof The idea of the proof is the same as in the original paper \cite{FeJiNo}. We will combine the energy inequality \eqref{eq:EI} provided by Proposition \ref{p:ex2} together with weak formulations of the continuity and momentum equations with suitable test functions. Since $\vU$ is not a proper test function in the momentum equation due to its boundary condition, we test the momentum equation with $\vph = \vU-\vV$ to obtain
\begin{align}\label{eq:pr1}
&\int_{\Omega_\tau} \vr \vu \cdot (\vU-\vV) (\tau, \cdot) \ \dx - \int_{\Omega_0} (\vr \vu)_0 \cdot (\vU-\vV) (0, \cdot) \ \dx \\ \nonumber
= &\int_0^\tau \int_{\Omega_t} \left( \vr \vu \cdot \partial_t (\vU-\vV) + \vr [\vu \otimes \vu] : \Grad (\vU-\vV) + p(\vr) \Div (\vU-\vV)
- \tn{S} (\Grad \vu) : \Grad (\vU-\vV)\right)\ \dxdt.
\end{align}
Subtracting \eqref{eq:pr1} from the energy inequality \eqref{eq:EI} we obtain
\begin{align}\label{eq:pr2}
&\int_{\Omega_\tau} \left(\frac 12\vr\abs{\vu}^2 + H(\vr) - \vr \vu \cdot \vU\right) (\tau, \cdot) \ \dx - \int_{\Omega_0} \frac{1}{2 \vr_0} |(\vr \vu)_0 |^2 + H(\vr_0) - (\vr \vu)_0 \cdot \vU (0, \cdot) \ \dx \\ \nonumber
+ &\int_0^\tau \int_{\Omega_t} \tn{S}(\Grad \vu):(\Grad\vu-\Grad\vU) \ \dxdt  \leq \int_0^\tau \int_{\Omega_t} \left( - \vr \vu \cdot \partial_t \vU - \vr [\vu \otimes \vu] : \Grad \vU - p(\vr) \Div \vU \right)\ \dxdt.
\end{align}
Next, we use in the continuity equation as a test function the quantities $\frac 12 \abs{\vU}^2$ and $H'(r)$ respectively to obtain
\bFormula{eq:pr3}
\int_{\Omega_\tau} \frac 12 \vr \abs{\vU}^2 (\tau, \cdot) \ \dx - \int_{\Omega_0} \frac 12 \vr_0 \abs{\vU}^2 (0, \cdot) \ \dx =
\int_0^\tau \int_{ \Omega_t} \left( \vr \vU \cdot \partial_t \vU + \vr \vu \cdot \Grad \vU \cdot \vU \right) \ \dxdt
\eF
and
\bFormula{eq:pr4}
\int_{\Omega_\tau} \vr H'(r) (\tau, \cdot) \ \dx - \int_{\Omega_0} \vr_0 H'(r) (0, \cdot) \ \dx =
\int_0^\tau \int_{ \Omega_t} \left( \vr \partial_t H'(r) + \vr \vu \cdot \Grad H'(r) \right) \ \dxdt.
\eF
Adding \eqref{eq:pr3} and subtracting \eqref{eq:pr4} from \eqref{eq:pr2} we obtain
\begin{align}\label{eq:pr5}
&\int_{\Omega_\tau} \left(\frac 12\vr\abs{\vu-\vU}^2 + H(\vr) - H'(r)\vr\right) (\tau, \cdot) \ \dx - \int_{\Omega_0} \frac{1}{2 \vr_0} |(\vr \vu)_0 - \vr_0\vU(0,\cdot) |^2 + H(\vr_0) - H'(r(0,\cdot))\vr_0 \ \dx \\ \nonumber
+ &\int_0^\tau \int_{\Omega_t} \tn{S}(\Grad \vu):(\Grad\vu-\Grad\vU) \ \dxdt  \leq \int_0^\tau \int_{\Omega_t} \left( (\vr \partial_t \vU + \vr\vu\cdot\Grad\vU)\cdot (\vU-\vu) - p(\vr) \Div \vU\right)\ \dxdt \\ \nonumber
- &\int_0^\tau \int_{\Omega_t} \left( \vr \partial_t H'(r) + \vr \vu \cdot \Grad H'(r) \right)\ \dxdt.
\end{align}
Observing that the definition \eqref{p1pb} implies
\begin{equation}\label{eq:Hpid}
p(r) = rH'(r) - H(r),
\end{equation}
we immediately achieve
\begin{equation}\label{eq:Hpid2}
\partial_t p(r) = r\partial_t H'(r).
\end{equation}
Hence, the inequality \eqref{eq:pr5} can be further rewritten as
\begin{align}\label{eq:pr6}
&\int_{\Omega_\tau} \left(\frac 12\vr\abs{\vu-\vU}^2 + H(\vr) - H'(r)\vr\right) (\tau, \cdot) \ \dx - \int_{\Omega_0} \frac{1}{2 \vr_0} |(\vr \vu)_0 - \vr_0\vU(0,\cdot) |^2 + H(\vr_0) - H'(r(0,\cdot))\vr_0 \ \dx \\ \nonumber
+ &\int_0^\tau \int_{\Omega_t} (\tn{S}(\Grad \vu) - \tn{S}(\Grad\vU)):(\Grad\vu-\Grad\vU) \ \dxdt + \int_0^\tau \int_{\Omega_t} \partial_t p(r)\ \dxdt  \\ \nonumber
\leq &\int_0^\tau \int_{\Omega_t} \left( (\vr \partial_t \vU + \vr\vu\cdot\Grad\vU)\cdot (\vU-\vu) + \tn{S}(\Grad\vU):(\Grad\vU-\Grad\vu)\right)\ \dxdt \\ \nonumber
+ &\int_0^\tau \int_{\Omega_t} \left((r - \vr) \partial_t H'(r) - p(\vr) \Div \vU - \vr \vu \cdot \Grad H'(r) \right)\ \dxdt.
\end{align}
Now we claim that the following identity holds
\begin{equation}\label{eq:pr7}
\int_{\Omega_t} p(r)\Div \vU  + r\vU\cdot\Grad H'(r) \ \dx = \int_{\Omega_t} \Div (\vV p(r)) \ \dx.
\end{equation}
Indeed, using the boundary condition $\vU\cdot\n = \vV\cdot\n$ we write
\begin{align}\label{eq:pr8}
\int_{\Omega_t} p(r)\Div\vU \ \dx &= \int_{\Omega_t} p(r)\Div(\vU-\vV)\ \dx + \int_{\Omega_t} p(r)\Div \vV \ \dx \\ \nonumber
&= -\int_{\Omega_t} \vU\cdot\Grad p(r)\ \dx + \int_{\Omega_t}\Div(\vV p(r))\ \dx = -\int_{\Omega_t} r\vU\cdot\Grad H'(r)\ \dx + \int_{\Omega_t}\Div(\vV p(r))\ \dx,
\end{align}
where we used \eqref{eq:Hpid} as well. Adding \eqref{eq:pr7} to \eqref{eq:pr6} we obtain
\begin{align}\label{eq:pr9}
&\int_{\Omega_\tau} \left(\frac 12\vr\abs{\vu-\vU}^2 + H(\vr) - H'(r)\vr\right) (\tau, \cdot) \ \dx - \int_{\Omega_0} \frac{1}{2 \vr_0} |(\vr \vu)_0 - \vr_0\vU(0,\cdot) |^2 + H(\vr_0) - H'(r(0,\cdot))\vr_0 \ \dx \\ \nonumber
+ &\int_0^\tau \int_{\Omega_t} (\tn{S}(\Grad \vu) -{ \tn{S}(\Grad\vu)):(\Grad\vU-\Grad\vU)} \ \dxdt + \int_0^\tau \int_{\Omega_t} (\partial_t p(r) + \Div (\vV p(r)))\ \dxdt  \\ \nonumber
\leq &\int_0^\tau \int_{\Omega_t} \left( (\vr \partial_t \vU + \vr\vu\cdot\Grad\vU)\cdot (\vU-\vu) + \tn{S}(\Grad\vU):(\Grad\vU-\Grad\vu)\right)\ \dxdt \\ \nonumber
+ &\int_0^\tau \int_{\Omega_t} \left((r - \vr) \partial_t H'(r) + (p(r)- p(\vr)) \Div \vU +(r\vU - \vr \vu) \cdot \Grad H'(r) \right)\ \dxdt.
\end{align}
The proof of Proposition \ref{p:REI} is finished observing that standard transport theorem yields the identity
\begin{align}\label{eq:pr10}
&\int_0^\tau \int_{\Omega_t} (\partial_t p(r) + \Div (\vV p(r)))\ \dxdt = \int_0^\tau \frac{{\rm d}}{\dt}\int_{\Omega_t} p(r)\ \dxdt \\ \nonumber
= &\int_{\Omega_\tau} p(r)(\tau,\cdot) \ \dx - \int_{\Omega_0} p(r)(0,\cdot) \ \dx = \int_{\Omega_\tau} (rH'(r)-H(r))(\tau,\cdot) \ \dx - \int_{\Omega_0} (rH'(r)-H(r))(0,\cdot) \ \dx.
\end{align}

\qed

Note that the class of admissible test functions $(r,\vU)$ can be extended by density arguments in a similar manner as in \cite[Section 3.2.2]{FeJiNo}.

\subsection{Proof of Theorem \ref{t:WSU}}\label{wsu}
The proof follows the same ideas as in \cite{FeJiNo}, however we present it here for completeness of presentation. Plugging in $(r,\vU) = (\vrt,\vut)$ in the relative energy inequality \eqref{eq:REI} we obtain
\begin{align}\label{eq:pr21}
&\Ecal\Big([\vr,\vu]|[\vrt,\vut]\Big)(\tau) + \int_0^\tau\int_{\Om_t} \left(\tn{S}(\Grad \vu) - \tn{S}(\Grad \vut)\right):\left(\Grad \vu - \Grad \vut\right)\,\dxdt \\ \nonumber &\leq \int_0^\tau \int_{\Om_t} \vr(\partial_t\vut + \vu\cdot\Grad\vut)\cdot(\vut-\vu) + \tn{S}(\Grad\vut):(\Grad\vut-\Grad\vu) \,\dxdt \\ \nonumber
&+ \int_0^\tau \int_{\Omega_t} \Div\vut(p(\vrt)-p(\vr)) + (\vrt-\vr)\partial_tH'(\vrt) + (\vrt\vut-\vr\vu)\cdot\Grad H'(\vrt)\, \dx \dt.
\end{align}
Using the strong formulation of the momentum and continuity equations we find out that
\begin{equation}\label{eq:pr22}
\partial_t\vut + \vut\cdot\Grad\vut = \frac{1}{\vrt} \Div \tn{S}(\Grad\vut) - \Grad H'(\vrt)
\end{equation}
in $Q_T$. Moreover, multiplying the strong formulation of the continuity equation by $H''(\vrt)$ we obtain
\begin{equation}\label{eq:pr23}
\partial_t H'(\vrt) + \vut\cdot\Grad H'(\vrt) = - \Div \vut \vrt H''(\vrt) = - \Div \vut p'(\vrt).
\end{equation}
Finally, integrating by parts we have for a.a. $t \in (0,\tau)$
\begin{equation}\label{eq:pr24}
\int_{\Omega_t} \tn{S}(\Grad\vut) : (\Grad\vut-\Grad\vu)\,\dx = - \int_{\Omega_t} \Div \tn{S}(\Grad\vut) \cdot (\vut-\vu)\,\dx,
\end{equation}
where the boundary integral vanishes due to boundary condition \eqref{i6c} (recall we set $\kappa = 0$) and the fact, that $(\vut-\vu)\cdot\n = 0$ on $\Gamma_t$. Combining \eqref{eq:pr21}, \eqref{eq:pr22}, \eqref{eq:pr23} and \eqref{eq:pr24} we arrive to the following version of the relative energy inequality
\begin{align}\label{eq:pr25}
&\Ecal\Big([\vr,\vu]|[\vrt,\vut]\Big)(\tau) + \int_0^\tau\int_{\Om_t} \left(\tn{S}(\Grad \vu) - \tn{S}(\Grad \vut)\right):\left(\Grad \vu - \Grad \vut\right)\,\dxdt \\ \nonumber &\leq \int_0^\tau \int_{\Om_t} \vr(\vu-\vut)\cdot\Grad\vut\cdot(\vut-\vu) - \Div\vut(p(\vr) - p'(\vrt)(\vr-\vrt) - p(\vrt))\,\dxdt \\ \nonumber
&+ \int_0^\tau \int_{\Omega_t} \frac{1}{\vrt}(\vr-\vrt)\Div\tn{S}(\Grad\vut)\cdot(\vut-\vu)\, \dx \dt.
\end{align}
Now we would like to show that all the terms on the right-hand side of \eqref{eq:pr25} can be absorbed by the left hand side and then use the Gronwall lemma. To do that we need the following estimate which can be easily checked
\begin{align}\label{eq:pr26}
H(\vr) - H'(r)(\vr-r) - H(r) &\geq c(r)(\vr-r)^2 \quad \text{ for } \frac{r}{2} < \vr < 2r \\ \nonumber
&\geq c(r)(1+\vr^\gamma) \quad \text{ otherwise}
\end{align}
and also the following Korn-type inequality
\begin{equation}\label{eq:Korn}
\norma{\vz}{W^{1,2}(\Omega_t)} \leq C\norma{\tn{S}(\Grad\vz)}{L^2(\Omega_t)}
\end{equation}
for all $\vz \in W^{1,2}(\Omega_t)$.

Thus, it is not difficult to observe that
\begin{equation}\label{eq:pr27}
\abs{\int_{\Om_t} \vr(\vu-\vut)\cdot\Grad\vut\cdot(\vut-\vu) - \Div\vut(p(\vr) - p'(\vrt)(\vr-\vrt) - p(\vrt))\,\dx} 
\leq C\norma{\Grad\vut}{L^\infty(\Omega_t)}\Ecal\Big([\vr,\vu]|[\vrt,\vut]\Big)(t).
\end{equation}
It remains to handle the last term on the right hand side of \eqref{eq:pr25}. We split the integral into three parts, considering first $\vr$ close to $\vrt$, then $\vr$ small and finally $\vr$ large. We have using the H\"older inequality, the Young inequality and \eqref{eq:Korn}
\begin{align}\label{eq:pr28}
&\abs{\int_{\{\vrt/2 \leq \vr \leq 2\vrt\}} \frac{1}{\vrt}(\vr-\vrt)\Div\tn{S}(\Grad\vut)\cdot(\vut-\vu)\,\dx} \\ \nonumber
\leq &C(\delta)\norma{\frac{1}{\vrt}\Div\tn{S}(\Grad\vut)}{L^3(\Omega_t)}^2\int_{\{\vrt/2 \leq \vr \leq 2\vrt\}}(\vr-\vrt)^2\,\dx + \delta\norma{\vut-\vu}{L^6(\Omega_t)}^2 \\ \nonumber
\leq &C(\delta)\norma{\frac{1}{\vrt}\Div\tn{S}(\Grad\vut)}{L^3(\Omega_t)}^2\Ecal\Big([\vr,\vu]|[\vrt,\vut]\Big)(t) + \delta C \norma{\tn{S}(\Grad(\vut)-\Grad(\vu))}{L^2(\Omega_t)}^2.
\end{align}
The last term can be absorbed into the left hand side for $\delta$ small enough, whereas the first term can be treated using the Gronwall lemma.

On the set where $\vr$ is small we can proceed in the following way
\begin{align}\label{eq:pr29}
&\abs{\int_{\{0\leq\vr<\vrt/2\}} \frac{1}{\vrt}(\vr-\vrt)\Div\tn{S}(\Grad\vut)\cdot(\vut-\vu)\,\dx}
\leq \abs{\int_{\{0\leq\vr<\vrt/2\}} \Div\tn{S}(\Grad\vut)\cdot(\vut-\vu)\,\dx} \\ \nonumber
&\qquad \leq C(\delta)\norma{\Div\tn{S}(\Grad\vut)}{L^3(\Omega_t)}^2\int_{\{0\leq\vr<\vrt/2\}}1\,\dx + \delta\norma{\vut-\vu}{L^6(\Omega_t)}^2 \\ \nonumber
&\qquad \leq C(\delta)\norma{\Div\tn{S}(\Grad\vut)}{L^3(\Omega_t)}^2\Ecal\Big([\vr,\vu]|[\vrt,\vut]\Big)(t) + \delta C \norma{\tn{S}(\Grad(\vut)-\Grad(\vu))}{L^2(\Omega_t)}^2.
\end{align}
Again, the last term is absorbed into the left hand side for $\delta$ small enough and the first term is treated using the Gronwall lemma.

Finally, consider the integral over the set where $\vr$ is large. Here we distinguish two cases. First for $\gamma \leq 2$ we have
\begin{align}\label{eq:pr30}
&\abs{\int_{\{\vr>2\vrt\}} \frac{1}{\vrt}(\vr-\vrt)\Div\tn{S}(\Grad\vut)\cdot(\vut-\vu)\,\dx} \leq \abs{\int_{\{\vr>2\vrt\}}\vr\frac{\vr-\vrt}{\vr\vrt}\Div\tn{S}(\Grad\vut)\cdot(\vut-\vu)\,\dx} \\ \nonumber
&\qquad \leq \abs{\int_{\{\vr>2\vrt\}}\vr\frac{1}{\vrt}\Div\tn{S}(\Grad\vut)\cdot(\vut-\vu)\,\dx} \\ \nonumber
&\qquad \leq C(\delta)\norma{\frac{1}{\vrt}\Div\tn{S}(\Grad\vut)}{L^{\frac{6\gamma}{5\gamma-6}}(\Omega_t)}^2\left(\int_{\{\vr>2\vrt\}}\vr^\gamma\,\dx\right)^{2/\gamma} + \delta\norma{\vut-\vu}{L^6(\Omega_t)}^2 \\ \nonumber
&\qquad \leq C(\delta)\norma{\frac{1}{\vrt}\Div\tn{S}(\Grad\vut)}{L^{\frac{6\gamma}{5\gamma-6}}(\Omega_t)}^2\Ecal\Big([\vr,\vu]|[\vrt,\vut]\Big)^{\frac{2}{\gamma}-1}(t)\Ecal\Big([\vr,\vu]|[\vrt,\vut]\Big)(t) + \delta C \norma{\tn{S}(\Grad(\vut)-\Grad(\vu))}{L^2(\Omega_t)}^2.
\end{align}
In this case the power $\frac{2}{\gamma}-1$ is nonnegative and we use also the property $\Ecal\Big([\vr,\vu]|[\vrt,\vut]\Big) \in L^\infty(0,T)$ to proceed further.

For $\gamma > 2$ we have
\begin{align}\label{eq:pr31}
&\abs{\int_{\{\vr>2\vrt\}} \frac{1}{\vrt}(\vr-\vrt)\Div\tn{S}(\Grad\vut)\cdot(\vut-\vu)\,\dx} \leq \abs{\int_{\{\vr>2\vrt\}}\vr\frac{\vr-\vrt}{\vr\vrt}\Div\tn{S}(\Grad\vut)\cdot(\vut-\vu)\,\dx} \\ \nonumber
&\qquad\leq \abs{\int_{\{\vr>2\vrt\}}\vr^{\frac{\gamma}{2}}\frac{1}{\vrt}\Div\tn{S}(\Grad\vut)\cdot(\vut-\vu)\,\dx} \\ \nonumber
&\qquad \leq C(\delta)\norma{\frac{1}{\vrt}\Div\tn{S}(\Grad\vut)}{L^3(\Omega_t)}^2\left(\int_{\{\vr>2\vrt\}}\vr^\gamma\,\dx\right) + \delta\norma{\vut-\vu}{L^6(\Omega_t)}^2 \\ \nonumber
&\qquad \leq C(\delta)\norma{\frac{1}{\vrt}\Div\tn{S}(\Grad\vut)}{L^3(\Omega_t)}^2\Ecal\Big([\vr,\vu]|[\vrt,\vut]\Big)(t) + \delta C \norma{\tn{S}(\Grad(\vut)-\Grad(\vu))}{L^2(\Omega_t)}^2.
\end{align}

Altogether we end up with the inequality
\begin{equation}\label{eq:Gron}
\Ecal\Big([\vr,\vu]|[\vrt,\vut]\Big)(\tau) \leq \int_0^\tau h(t)\Ecal\Big([\vr,\vu]|[\vrt,\vut]\Big)(t)\,\dt
\end{equation}
for some $h(t) \in L^1(0,T)$ and the Gronwall lemma finishes the proof.

\qed

\begin{Remark}
The case of nonzero bulk viscosity coefficient $\eta > 0$ in \eqref{i4} as well as the case of nonzero boundary friction coefficient $\kappa > 0$ in \eqref{i6c} can be treated by obvious modifications just adding proper integrals to appropriate formulas.
\end{Remark}

\section{Concluding remarks}\label{s:conc}
For clarity of the proof we have so far restricted our presentation to the case of slip boundary conditions. In case of homogeneous Dirichlet boundary condition
\bFormula{i6a}
(\vu - \vc{V})|_{\Gamma_\tau} = 0 \ \mbox{for any}\ \tau \geq 0
\eF
the weak-strong uniqueness principle has been shown recently in \cite{dobo}. However, the existence of regular solutions has remained so far open question. Theorem \ref{Tmain} easily extends to this case as well, in fact the proof can be not only repeated but considerably simplified since the boundary condition \eqref{i6a} remains homogeneous under Lagrangian transformation and therefore we do not need the extension operator defined in Lemma \ref{Lext}.

It is also clear from our proof that it remains valid if we assume the right hand side in the momentum equation on the form $\vr \vc{f}$ with
$$
\vc{f} \in L^2(0,T;H^1(\Omega_t)), \quad \vc{f}_t \in L^2(0,T;L^2(\Omega_t)), \quad \vc{f}|_{t=0} \in H^1(\Omega).
$$

We also notice that the regularity assumptions on $\vV$ in Theorem \ref{Tmain} are not optimal. However, we need some integrability of the third order derivatives of $\vV$ and therefore it is not enough to assume the regularity from Theorem \ref{t:WSU} which is sufficient for the existence of weak solutions and for the weak-strong uniqueness.

Taking into account known existence results for weak solutions \cite{FeNeSt}, \cite{FKNNS} and the weak-strong uniqueness result \cite{dobo}, our paper completes a part of the local existence theory for the compressible barotropic Navier-Stokes system on moving domains at least in the framework of Hilbert spaces. A natural generalization now could be existence result for regular solutions in $L^p$ setting. A more involved interesting issue is the global well-posedness of strong solutions for small data. In case of the complete system with thermal effects, for which existence of weak solutions on moving domains has been shown recently in \cite{KMNW1} and \cite{KMNW2}, both the existence of strong solutions and weak-strong uniqueness remain open problems.

\appendix
\begin{appendices}
\renewcommand{\theequation}{\Alph{section}.\arabic{equation}}
\renewcommand\thesection{}
%
\section{Appendix}
{\bf Proof of Lemma \ref{Lext}.}
For the purpose of our construction it is convenient to define the whole
velocity at the boundary. Therefore we look for the extension of the boundary data satisfying conditions
\begin{align} \label{ext:10}
&\vu^b(t,y) \cdot \vc{n}(y) = (\vu - \vc{V})(t,y)\cdot\left(\vc{n}(y)-\vc{n}(\vc{X}(t,y))\right) + (\vc{V}(t,\vc{X}(t,y))-\vc{V}(t,y))\cdot \vc{n}(\vc{X}(t,y)),\\ \nonumber
&\vu^b(t,y) \cdot \tau^k(y) = (\vu - \vc{V})(t,y)\cdot\left(\tau^k(y)-\tau^k(\vc{X}(t,y))\right) + (\vc{V}(t,\vc{X}(t,y))-\vc{V}(t,y))\cdot \tau^k(\vc{X}(t,y))
\end{align}
and
\begin{align} \label{ext:10b}
&\mu(\nabla_y \vu^b+\nabla^T_y \vu^b)(t,y)\vc{n}(y) \cdot \tau^k(y)=
\mu\left(\nabla_y\vu(t,y)(\tn{I}-\nabla_x\vc{Y}) + ((\tn{I}-\nabla_x^T\vc{Y})\nabla_y^T\vu(t,y))^T\right)\vc{n}(\vc{X}(t,y))\cdot\tau^k(\vc{X}(t,y)) \\ \nonumber
& \qquad + \mu(\nabla_y\vu+\nabla^T_y\vu)(t,y)[(\vc{n}(y)-\vc{n}(\vc{X}(t,y)))\cdot\tau^k(\vc{X}(t,y)) + \vc{n}(t,y)\cdot(\tau^k(y)-\tau^k(\vc{X}(t,y)))].
\end{align}
First of all, notice that it is enough to define appropriate extension only
in a neighbourhood of the boundary
$$
\Omega_{\epsilon}=\{x \in \Omega: {\rm dist}(x,\de \Omega)<\epsilon\}.
$$
Then multiplying it by a smooth
function $\phi$ such that
$$
\phi(x)\in[0,1], \quad \phi|_{\Omega_{\epsilon}}\equiv 1, \quad \phi_{\Omega \setminus \Omega_{2\epsilon}}\equiv 0
$$
we obtain a function defined on the whole $\Omega$ which also satisfies
the estimate \eqref{est_ext}. Next important observation is that it is enough
to consider the case of flat boundary which is obtained by nowadays classical technique of partition of unity. The cutoff functions involved in this procedure enjoys the regularity of the boundary of $\Omega$ and therefore all their contribution to our estimates can be put in the constant in \eqref{est_ext}. Therefore we assume
\begin{equation} \label{ext:7}
\vn(x_1,x_2,0)=(0,0,1), \quad \tau^1(x_1,x_2,0)=(1,0,0), \quad \tau^2(x_1,x_2,0)=(0,1,0).
\end{equation}
Our construction will be carried out in two steps. First we find $\tilde \vu^b$
satisfying only relations \eqref{ext:10}. In the second step we will use it
to define an extension satisfying also the relation for the derivatives.
Let us focus on the extension of the normal component of $\tilde \vu^b$.
For $y=(y_1,y_2,y_3)\in \Omega$  we define
$$
\delta \vn(t,y)=\vn(t,(y_1,y_2,0))-\vn(\vX(t,(y_1,y_2,0)))
$$
and
$$
\delta \vV(t,y)=\vV(t,\vX(t,(y_1,y_2,0)))-\vV(t,(y_1,y_2,0)).
$$
Then it natural do define the extension of the normal component of $\vu^b$ as
\begin{equation} \label{ext:0}
\vu^b_1(t,y)=(\vu-\vV)(t,y)\cdot \delta \vn(t,y) + \delta \vV(t,y)\cdot \vn(\vX(t,(y_1,y_2,0))).
\end{equation}
Let us start with the first component. We have
\begin{equation} \label{ext:1}
[(\vu-\vV)\cdot \delta \vn]_{tt}=(\vu-\vV)_{tt}\delta \vn + 2(\vu-\vV)_t(\delta \vn)_t+(\vu-\vV)(\delta \vn)_{tt}.
\end{equation}
By \eqref{eq:coc} we have
$$
\delta \vn(t,y) \sim \int_{0}^t \vV, \quad (\delta \vn(t,y))_t \sim \vV, \quad (\delta \vn(t,y))_{tt}\sim \vV_t.
$$
We use these relations to estimate the $L^2(L^2)$ norm of the right hand side of \eqref{ext:1}. The idea is that when we have $(\vu-\vV)_{tt}$
we can get the smallness in time from $\int_0^t \vV$ and in the remaining terms we get smallness in time using boundedness in time of appropriate norms of $\vu-\vV$. Precisely,
for the first term we have
$$
\int_0^T\|(\vu-\vV)_{tt}\delta \vn\|_{L^2(L^2)}^2\leq \int_0^T\|\delta\vn\|_{L^\infty}^2 \|(\vu-\vV)_{tt}\|_{L^2}^2 \leq T^2 \|\vV\|_{L^\infty(\Omega \times (0,T))}^2\int_0^T\|(\vu-\vV)_{tt}\|_{L^2}^2,
$$
for the second
$$
\int_0^T \|(\vu-\vV)_t(\delta \vn)_t\|_{L^2}^2 \leq \|(\delta\vn)_t\|_{L^\infty(\Omega \times (0,T))}^2 \int_0^T \|(\vu-\vV)_t\|_{L^2}^2 \leq T \|\vV\|_{L^\infty(\Omega \times (0,T))}\|(\vu-\vV)_t\|_{L^\infty(L^2)}^2
$$
and similarly for the third with $\|\vV_t\|_{L_\infty}$.
We conclude
\begin{equation} \label{ext:2}
\|[(\vu-\vV)\delta \vn]_{tt}\|_{L^2(L^2)} \leq C(T+\sqrt{T})\|\vV\|_{W^1_\infty(L_\infty)}[\|(\vu-\vV)_{tt}\|_{L^2(L^2)} +\|(\vu-\vV)\|_{W^1_\infty(L^2)}].
\end{equation}
Now we write the second time derivative of the second term in \eqref{ext:0} (we denote $E\vn:=\vn(\vX(t,(y_1,y_2,0))$:
\begin{equation} \label{ext:3}
(\delta \vV \cdot E\vn)_{tt}=(\delta \vV)_{tt}\cdot E\vn + 2 (\delta \vV)_t (E\vn)_t + (E\vn)_{tt}.
\end{equation}
Notice that we have
$$
\delta \vV \sim \nabla \vV (\vX(t,(y_1,y_2,0))-(y_1,y_2,0)) \sim \nabla \vV \int_0^t \vV.
$$
Therefore
$$
(\delta \vV)_t \sim \nabla \vV_t \int_0^t \vV + \vV \nabla \vV
$$
and
$$
(\delta \vV)_{tt} \sim \nabla \vV_{tt}\int_0^t \vV + \vV \nabla \vV_t + \vV_t \nabla \vV.
$$
Like previously, in the first term of the RHS of the latter formula we get the smallness in time from $\int_0^t \vV$:
$$
\int_0^T \|\nabla \vV_{tt} \int_0^t \vV E\vn\|_{L^2}^2 {\rm d}t \leq CT\|\vV\|_{L_\infty(\Omega \times (0,T))}\|\nabla \vV_{tt}\|_{L^2(L^2)}.
$$
In the terms where we does not have this factor we still get smallness in time using boundedness in time of appropriate norms of $\vV$, for example
$$
\int_0^t \|\vV \cdot \nabla \vV (E\vn)_t\|_{L^2}^2\leq Ct \|\vV\|^4_{L^\infty(W^1_\infty)}.
$$
Treating similarly the other terms in \eqref{ext:3} we obtain
\begin{equation}
\|\tilde \vu^b_{1,tt}\|_{L^2(L^2)} \leq C \sqrt{T}(1+\sqrt{T})\phi(\|\vV\|_{L^\infty(W^1_\infty)},\|\vV\|_{W^1_\infty(L_\infty)},\|\nabla \vV_{tt}\|_{L^2(L^2)}).
\end{equation}
Combining this estimate with \eqref{ext:2} we conclude
\begin{equation} \label{ext:4}
\|E(\vu^b \cdot \vn)_{tt}\|_{L^2(L^2)} \leq E(t)[1+\|(\vu-\vV)_{tt}\|_{L^2(L^2)}+\|\vu-\vV\|_{W^1_\infty(L^2)}].
\end{equation}
The extension of tangential components is done exactly in the same way and so the estimate \eqref{ext:4} holds for the whole $\tilde \vu^b$. 
Now we will use $\tilde \vu^b$ to construct a function with the same values on the boundary satisfying also the relations \eqref{ext:10b} for tangential stress.
Conditions \eqref{ext:10} can be written in a compact form
\begin{equation} \label{ext:6}
\vu^b = E^1 \vu + E^2
\end{equation}
where $E^1$ and $E^2$ are small and sufficiently regular matrix and vector functions respectively. Furthermore, conditions \eqref{ext:10b} rewrite as (we denote $f_{y_i}$ by $f_{,i}$):
\begin{align} \label{ext:8}
&u^b_{1,3}+u^b_{3,1}= \sum_{i,j=1}^3A_{ij}(t,x)u_{i,j}\\
&u^b_{2,3}+u^b_{3,2}= \sum_{i,j=1}^3B_{ij}(t,x)u_{i,j}.
\end{align}
First of all it is natural to take $u^b_3=\tilde u^b_3$.
Next we can construct $u^b_1$ and $u^b_2$ separately. As both will be defined  analogously, we focus on $u^b_1$. As $u^b$ is determined on the boundary, so are its tangential derivatives, hence $\de_{x_1}$ and $\de_{x_2}$ due to \eqref{ext:7}. In particular,
$$
u^b_{3,1}=\tilde u^b_{3,1}=\sum_{i=1}^3E^{1}_{3i}u_{i,1}+\sum_{i=1}^3 E^1_{3i,1}u_i+E^2_{3,1}.
$$
Substituting this relation to \eqref{ext:8} we get
\begin{equation} \label{ext:9}
u^b_{1,3}=\sum_{i,j=1}^3 \bar A_{ij}u_{i,j}-\sum_i E^1_{3i,1}u_i-E^2_{3,1}.
\end{equation}
Now we can define $u^b_1=u^{b1}_1+u^{b2}_1$ where
\begin{equation} \label{ext:11}
u^{b1}_1(y)= (\sum_{i=1}^3E^1_{1i} u_i + E^2_1)(y)+ 2 \sum_{i,j}\tilde A_{ij}u_i((y_1,y_2,0)+y_3 {\bf e}_j)-2 \sum_{i,j}\tilde A_{ij}u_i((y_1,y_2,0)+\frac{y_3}{2} {\bf e}_j),
\end{equation}
where
\begin{align} \label{ext:12}
&\tilde A_{ij}=\bar A_{ij}, \qquad j\neq3,\nonumber\\
&\tilde A_{i3}=\bar A_{i3} - E^1_{1,i}, \quad  i=1,2,3.
\end{align}
Differentiating \eqref{ext:11} wrt $x_3$ we obtain
$$
u^{b1}_{1,3}=\sum_i E^1_{1i} u_{i,3}+\sum_i E^1_{1i,3}u_i+E^2_{1,3}+\sum_{i,j}\tilde A_{ij} u_{i,j},
$$
which by the definition of $\tilde A_{ij}$ reduces to
\begin{equation} \label{ext:13}
u^{b1}_{1,3}=\sum_{i,j}\bar A_{ij} u_{i,j}+\sum_i E^1_{1i,3}u_i+E^2_{1,3}.
\end{equation}
We see that $u^{b1}=\tilde u^{b1}$ on the boundary and \eqref{ext:13}
differs from \eqref{ext:9} only up to lower order terms (without derivatives of $u$). Therefore we need $u^{b2}_1=0$ on the boundary which will compensate these lower order terms, hence it must satisfy
$$
u^{b2}_{1,3} = -\sum_i E^1_{1i,3}u_i-E^2_{1,3}-\sum_i E^1_{3i,1}u_i-E^{2}_{3,1}=:P_u(x).
$$
We can define $u_1^{b2}$ simply as
\begin{equation} \label{ext:14}
u^{b2}_{1,3} = \int_0^{x_3} P_u(x_1,x_2,s)ds.
\end{equation}
We see that $u^{b1}=u^{b1}_1+u^{b1}_2$ satisfies all required relations.
The estimate \eqref{ext:4}
for $u^{b1}_1$ can be shown similarly to the estimate for $\tilde u^b$.
Namely, the structure of the terms $\tilde A_{ij}$ is either of a form $\tau(y)-\tau(\vc{X}(y))$ which has been treated in the estimate for $\tilde u^b$ or of a form $\tn{I}-\nabla_x\vc{Y}$ which can be treated in the same may. Indeed, we observe easily that
$$
\|\tn{I}-\nabla_x\vc{Y}\|_{L^\infty(\Omega \times (0,T))}\leq E(T)
$$
where $E$ is continuous, $E(0)=0$ and
$$
\|(\tn{I}-\nabla_x\vc{Y})_t,(\tn{I}-\nabla_x\vc{Y})_{tt}\|_{L^\infty(\Omega \times (0,T))}\leq C.
$$
Therefore, we obtain the estimate for $u^{b1}_1$. In order to show the estimate for $u^{b2}_1$ it is enough to show it for the integrand $P_u$ due to boundedness of $\Omega$. Again, $P_u$ contains the terms of structure which we already investigated except for space derivatives of $E^1$. However, we have
$$
\nabla E^1,\nabla E^2 \sim \int_0^t \nabla \vV,
$$
therefore we can repeat the estimates having the assumed regularity of $\vV$.
Finally, the second component of $u^b$ is constructed in the same way. Therefore all the arguments can be repeated and we obtain the estimate
for $\|\vu^b_{tt}\|_{L^2(L^2)}$. As now we have the explicit formula for $\vu^b$, the other norms entering $\|\vu^b_{tt}\|_{{\cal{Y}}(T)}$ are obtained using similar arguments, therefore we skip the details.

\qed

{\bf Penalized problem - weak formulation.}
For convenience of the reader we recall the approximation scheme from work \cite{FKNNS}.

The system of equations is considered in a fixed spatial domain $B$ which is large enough so that $\Omega_t \subset B$ for all $t \in [0,T)$. Then the time-space domain $(0,T)\times B$ is decomposed into the fluid part $Q_T$ and the solid part $Q_T^c := ((0,T)\times B) \setminus Q_T$. The authors add to the variational (weak) formulation of the momentum equation  a singular forcing term
\bFormula{i8new}
\frac{1}{\ep} \int_0^T \int_{\Gamma_t}  (\vu - \vc{V} ) \cdot \vc{n} \ \vph \cdot \vc{n} \ {\rm dS}_{x} \ \dt,\ \ep > 0 \ \mbox{small},
\eF
which penalizes the normal component of the velocity on the boundary of the fluid domain.

In \cite{FKNNS} the following three level penalization scheme was introduced

\begin{enumerate}
\item In addition to (\ref{i8new}), the authors introduced a \emph{variable} shear viscosity coefficient $\mu = \mu_\omega$, where $\mu_\omega$ remains strictly positive in the fluid domain $Q_T$ but vanishes in the solid domain $Q_T^c$ as $\omega \to 0$.

\item Similarly to the existence theory developed in \cite{FNP}, the authors introduced the \emph{artificial pressure}
\[
p_\delta(\varrho) = p(\varrho) + \delta \varrho^\beta,\ \beta \geq 2,\ \delta > 0,
\]
in the momentum equation (\ref{i1b}).
\item The initial data $\vr_0$ and $(\vr\vu)_0$ are prolonged from $\Omega_0$ to $B$ by zero and in $\Omega_0$ regularized to $\vr_{0,\delta}$, $(\vr\vu)_{0,\delta}$ such that  $\int_{B}
 \left( \vr_{0, \delta}^\gamma +
 \delta \vr_{0, \delta}^\beta \right) \dx \leq c$ and $\int_{\Omega_0} \frac{1}{\vr_{0, \delta}} |(\vr \vu)_{0,\delta} |^2 \ \dx \leq c$.
\item Keeping $\ep, \delta,\ \omega > 0$ fixed, the authors solved the modified problem in $(0,T)\times B$. More precisely the weak solution to the penalized problem consists of couple $(\vr,\vu)$ satisfying
\bFormula{p3new}
\int_{B} \vr \varphi (\tau, \cdot) \ \dx - \int_{B} \vr_{0,\delta} \varphi (0, \cdot) \ \dx =
\int_0^\tau \int_{B} \left( \vr \partial_t \varphi + \vr \vu \cdot \Grad \varphi \right) \ \dxdt
\eF
for any $\tau \in [0,T]$ and any test function $\varphi \in \DC([0,T] \times \R^3)$;
\bFormula{p4new}
\int_{B} \vr \vu \cdot \vph (\tau, \cdot) \ \dx - \int_{B} (\vr \vu)_{0,\delta} \cdot \vph (0, \cdot) \ \dx
\eF
\[
= \int_0^\tau \int_{B} \left( \vr \vu \cdot \partial_t \vph + \vr [\vu \otimes \vu] : \Grad \vph + p(\vr) \Div \vph + \delta \vr^\beta \Div \vph
- \mu_\omega \left( \Grad \vu + \Grad^t \vu - \frac{2}{3} \Div \vu \tn{I} \right) : \Grad \vph \right) \ \dxdt
\]
\[
+ \frac{1}{\ep} \int_0^\tau \int_{ \Gamma_t } \left( (\vc{V} - \vu ) \cdot \vc{n} \ \vph \cdot \vc{n} \right) \ {\rm dS}_{x} \ \dt
\]
for any $\tau \in [0,T]$ and any test function $\vph \in \DC([0,T] \times B ; \R^3)$,
where $\vu \in L^2(0,T; W^{1,2}_0 (B; \R^3))$, meaning
$\vu$ satisfies the no-slip boundary condition
\bFormula{p5new}
\vu|_{\partial B} = 0 \ \mbox{in the sense of traces}.
\eF

\item Finally, the authors performed successively the limits $\varepsilon \rightarrow 0$, $\omega \rightarrow 0$ and $\delta \rightarrow 0$ to obtain a weak solution as in Definiton \ref{d:ws}.

\end{enumerate}
\end{appendices}

\smallskip
\noindent
{\bf Acknowledgements.} The authors would like to thank Milan Pokorn\'y for valuable suggestions concerning the proof of Lemma \ref{Lext}.
Moreover, they would like to thank to anonymous referee for her/his comments that helped us to improve the maniscript.

\end{document}